%% file: paper4.tex

\input makro.tex

\input smalmtrx.tex

\duplex

\def\star{{\relax}\tilde{\ }}

\title{MORITA EQUIVALENCE OF $C^{*}$-CROSSED PRODUCTS}
\title{BY INVERSE SEMIGROUP ACTIONS AND PARTIAL ACTIONS}

\author{N\'andor Sieben}
\shorttitle{morita equivalent semigroup actions}

\footnote{}{
\vskip-\baselineskip

\noindent This material is based upon work supported
by the National Science Foundation under Grant No.
DMS9401253.}

\bigskip\abstract{Morita equivalence of twisted inverse semigroup 
actions and discrete twisted partial actions are introduced.  Morita equivalent 
actions have Morita equivalent crossed products.  } 

\newsection{Introduction}

Morita equivalence of group actions on $C^{*}$-algebras was studied by 
Combes [Com], Echterhoff [Ech], Curto, Muhly and Williams [CMW] and 
Kaliszewski [Kal].  We adapt this notion for both Busby-Smith and 
Green type inverse semigroup actions, introduced in [Si1] and [Si2].  
We show that Morita equivalence is an equivalence relation and that 
Morita equivalent actions have Morita equivalent crossed products.  
The close connection between inverse semigroup actions and partial 
actions [Si1], [Ex3], [Si2] makes it easy to find the notion of Morita 
equivalence for discrete twisted partial actions.  In Section 4 we 
work out some of the details of discrete twisted partial crossed 
products, continuing the work started in [Ex2].  The fact that Morita 
equivalent twisted partial actions have Morita equivalent crossed 
products will then follow from the connection with semigroup 
actions.  In [AEE] Abadie, Eilers and Exel introduced Morita 
equivalence of crossed products by Hilbert bimodules.  We show that 
this definition is equivalent to our definition of Morita equivalence 
on the common special case of partial actions by {\bf Z}.  

The research for this paper was carried out while the author was a 
student at Arizona State University under the supervision of John 
Quigg.  I thank Professor Quigg for his help during the 
writing of this paper.

\newsection{Preliminaries}

\uj In this section we recall some basic definitions to fix our 
terminology and notation. Our references for Hilbert modules are 
[JT] and [Lan].  

Let $B$ be a $C^{*}$-algebra.  A ({\it right\/}) $B${\it -module\/} is a 
complex vector space $X$ with a bilinear map 
$(x,b)\mapsto x\cdot b:X\times B\to X$ such that $(x\cdot b)\cdot 
c=x\cdot (b\cdot c)$ for all 
$x\in X$ and $b,c\in B$.  
A ({\it right\/}) {\it inner\/}-{\it product }
$B${\it -module\/} is a $B$-module with a map 
$\langle\cdot ,\cdot\rangle_B:X\times X\to B$, called a $B${\it -valued inner product},
such that for all $\lambda ,\mu\in {\bf C}$, $x,y,z\in X$ and 
$b\in B$ we have
\item{(a)}{$\langle x,\lambda y+\mu z\rangle_B=\lambda\langle x,y
\rangle_B+\mu\langle x,z\rangle_B$;} 
\item{(b)}{$\langle x,y\cdot b\rangle_B=\langle x,y\rangle_Bb$;} 
\item{(c)}{$\langle x,y\rangle_B^{*}=\langle y,x\rangle_B$;} 
\item{(d)}{$\langle x,x\rangle_B\ge 0$;} 
\item{(e)}{$\langle x,x\rangle =0$ only if $x=0$.} 

In an inner product $B$-module we have a norm 
$\|x\|_B=\|\langle x,x\rangle_B\|^{1/2}$ satisfying $\|x\cdot b\|_
B\leq\|x\|_B\|b\|$ and 
$\|\langle x,y\rangle_B\|\le\|x\|\,\|y\|$ for all $x,y\in X$ and $
b\in B$.  
A ({\it right\/}) {\it Hilbert} $B${\it -module\/} is an inner-product 
$B$-module, which is complete in the norm $\|\cdot\|_B$.  A 
Hilbert $B$-module $X$ satisfying 
$$\hbox{\rm $\overline {\hbox{\rm span}}$}\{\langle x,y\rangle_B:
x,y\in X\}=B$$
is called {\it full}.

Left modules are defined similarly, with the left inner 
product linear in the first variable.  For a left inner-product 
$A$-module we use the notation $_A\langle\cdot ,\cdot\rangle$ for the $
A$-valued inner product.  

\mark{fullemma}
\lemma Let $X$ be a full right Hilbert $B$-module and $b\in B$. 
If $x\cdot b=0$ for all $x\in X$, then $b=0.$

\proof For all $x,y\in X$ we have $\langle x,y\rangle_Bb=\langle 
x,y\cdot b\rangle_B=0$ 
which implies that $b=0$ by fullness.  \eop

\uj
Let $A$ and $B$ be $C^{*}$-algebras. An 
$A-B$-{\it bimodule} $_AX_B$ is a right $B$-module $X$ which is also a left 
$A$-module satisfying 
$$a\cdot (x\cdot b)=(a\cdot x)\cdot b$$
for all $a\in A$, $x\in X$ and $b\in B$.  
Note that a bimodule satisfies $(\lambda a)\cdot (x\cdot b)=a\cdot 
(x\cdot (\lambda b))$ for all $\lambda\in {\bf C}$.

\definition Let $A$ and $B$ be $C^{*}$-algebras.  A {\it Hilbert }
$A-B${\it -bimodule\/} is a bimodule $_AX_B$ which is a left Hilbert 
$A$-module and a right Hilbert $B$-module such that 
 
$$_A\langle x,y\rangle\cdot z=x\cdot\langle y,z\rangle_B$$
for all $x,y,z\in X$.  

\noindent A Hilbert bimodule which is also full on both 
sides is called an {\it imprimitivity bimodule. }

\uj Note that for any Hilbert bimodule $_AX_B$ there is a 
corresponding imprimitivity bimodule $_{A_0}X_{B_0}$ where 
$A_0=\overline {\hbox{\rm span}}\,{}_A\langle X,X\rangle$ and $B_
0=\overline {\hbox{\rm span}}\,\langle X,X\rangle_B$.

\lemma If $_AX_B$ is a Hilbert bimodule then  
\item{{\rm(a)}}{$\langle a\cdot x,y\rangle_B=\langle x,a^{*}\cdot 
y\rangle_B${\rm ;}} 
\item{{\rm (b)}}{$_A\langle x\cdot b,y\rangle ={}_A\langle x,y\cdot b^{
*}\rangle ,$} 

\noindent
for all $a\in A$, $b\in B$ and $x,y\in X$.

\proof Part (a) follows from the following calculation:
$$\eqalign{&\|\langle a\cdot x,y\rangle_B-\langle x,a^{*}\cdot y\rangle_
B\|^2\cr
&\qquad =\|(\langle a\cdot x,y\rangle_B-\langle x,a^{*}\cdot y\rangle_
B)^{*}(\langle a\cdot x,y\rangle_B-\langle x,a^{*}\cdot y\rangle_
B)\|\cr
&\qquad =\|\langle y,a\cdot x\rangle_B\langle a\cdot x,y\rangle_B
+\langle a^{*}\cdot y,x\rangle_B\langle x,a^{*}\cdot y\rangle_B\cr
&\qquad\qquad\qquad -\langle y,a\cdot x\rangle_B\langle x,a^{*}\cdot 
y\rangle_B-\langle a^{*}\cdot y,x\rangle_B\langle a\cdot x,y\rangle_
B\|\cr
&\qquad =\|\langle y,a\cdot x\cdot\langle a\cdot x,y\rangle_B\rangle_
B+\langle a^{*}\cdot y,x\cdot\langle x,a^{*}\cdot y\rangle_B\rangle_
B\cr
&\qquad\qquad\qquad -\langle y,a\cdot x\cdot\langle x,a^{*}\cdot 
y\rangle_B\rangle_B-\langle a^{*}\cdot y,x\cdot\langle a\cdot x,y
\rangle_B\rangle_B\|\cr
&\qquad =\|\langle y,a_A\langle x,a\cdot x\rangle y\rangle_B+\langle 
a^{*}\cdot y,_A\langle x,x\rangle a^{*}\cdot y\rangle_B\cr
&\qquad\qquad\qquad -\langle y,a_A\langle x,x\rangle a^{*}\cdot y
\rangle_B-\langle a^{*}\cdot y,_A\langle x,a\cdot x\rangle\cdot y
\rangle_B\|=0\,.\cr}
$$
Part (b) can be proved similarly.
\eop

\mark{impiso} 
\definition The triple $(\phi_A,\phi ,\phi_B)$ is called an {\it isomorphism between the }
{\it Hilbert bimodules} $_AX_B$ and $_CY_D$ if $\phi_A:A\to C$ and $
\phi_B:B\to D$ are 
$C^{*}$-algebra isomorphisms and $\phi :X\to Y$ is a map such that for all 
$x,y\in X$ and $a\in A$, $b\in B$ we have 
\item{(a)}{$\phi (x\cdot b)=\phi (x)\cdot\phi_B(b)$;} 
\item{(b)}{$\phi_B(\langle x,y\rangle_B)=\langle\phi (x),\phi (y)
\rangle_D$;} 
\item{(c)}{$\phi (a\cdot x)=\phi_A(a)\cdot\phi (x)$;} 
\item{(d)}{$\phi_A(_A\langle x,y\rangle )={}_C\langle\phi (x),\phi 
(y)\rangle$;} 
\item{(e)}{$\phi$ is surjective.} 

\uj The following lemma shows that we can relax some of these 
conditions. Note that part (ii) is an improvement of [Kal, Lemma 
1.1.3].

\mark{relaxlem} 

\lemma With the notations of Definition \cite{impiso} we have 
\item{{\rm(i)}}{if $\phi$ satisfies {\rm(b)} then it is a linear isometry$
;$} 
\item{{\rm(ii)}}{if $\phi$ satisfies {\rm(b)} then it also satisfies 
{\rm(a);}} 
\item{{\rm(iii)}}{if $\phi$ satisfies {\rm(b)} and {\rm(c)} and $_
CY_D$ is an 
imprimitivity bimodule then $\phi$ also satisfies {\rm(d)} and {\rm(e)} so that it is an 
isomorphism between $X$ and $Y$.} 

\proof An easy calculation using (b) and the linearity of $\phi_B$ shows that 
$\|\phi (\lambda a+\mu b)-\lambda\phi (a)-\mu\phi (b)\|^2=0$. 
It is also an isometry since 
$$\eqalign{\|\phi (x)\|^2=\|\langle\phi (x),\phi (x)\rangle_D\|=\|
\phi_B(\langle x,x\rangle_B)\|=\|\langle x,x\rangle_B\|=\|x\|^2\,
.\cr}
$$
Part (ii) follows from the following calculation:
$$\eqalign{&\|\phi (x\cdot b)-\phi (x)\cdot\phi_B(b)\|^2\cr
&\qquad =\|\langle\phi (x\cdot b)-\phi (x)\cdot\phi_B(b),\phi (x\cdot 
b)-\phi (x)\cdot\phi_B(b)\rangle_D\|\cr
&\qquad =\|\langle\phi (x\cdot b),\phi (x\cdot b)\rangle_D-\langle
\phi (x\cdot b),\phi (x)\cdot\phi_B(b)\rangle_D\cr
&\qquad\qquad -\langle\phi (x)\cdot\phi_B(b),\phi (x\cdot b)\rangle_
D+\langle\phi (x)\cdot\phi_B(b),\phi (x)\cdot\phi_B(b)\rangle_D\|\cr
&\qquad =\|\phi_B(\langle x\cdot b,x\cdot b\rangle_B)-\phi_B(\langle 
x\cdot b,x\rangle_B)\phi_B(b)\cr
&\qquad\qquad -\phi_B(b^{*})\phi_B(\langle x,x\cdot b\rangle_B)+\phi_
B(b^{*})\phi_B(\langle x,x\rangle_B)\phi_B(b)\|=0\,.\cr}
$$
To show (iii) let $Z=\overline {\phi (X)}$. Then we have 
$$\eqalign{D=\phi_B(B)&=\phi_B(\overline {\hbox{\rm span}}\,\langle 
X,X\rangle_B)\subset\overline {\hbox{\rm span}}\,\phi_B(\langle X
,X\rangle_B)\cr
&=\overline {\hbox{\rm span}}\,\langle\phi (X),\phi (X)\rangle_D\subset\overline {\hbox{\rm span}}\,
\langle Z,Z\rangle_D\,,\cr}
$$
and so $D=\overline {\hbox{\rm span}}\,\langle Z,Z\rangle_D$.
$Z$ is a left $C$-module since
$$\eqalign{C\cdot Z&=\phi_A(A)\cdot\overline {\phi (X)}\subset\overline {
(\phi_A(A)\cdot\phi (X))}\cr
&=\overline {\phi (A\cdot X)}=\overline {\phi (X)}=Z\,.\cr}
$$
$Z$ is also a right $D$-module since
$$\eqalign{Z\cdot D&\subset\overline {\phi (X)}\cdot\overline {\hbox{\rm span}}\,
\langle\phi (X),\phi (X)\rangle_D\cr
&=\overline {\hbox{\rm span}}\,(_C\langle\phi (X),\phi_{}(X)\rangle
\cdot\phi (X))\subset Z\,.\cr}
$$
Hence $Z$ is a closed subbimodule of $Y$ with full right inner product, 
and so $Z=Y$ by the Rieffel correspondence.  This shows that 
$\phi (X)=Y$ since $\phi$ is an isometry between Banach spaces.  For 
$x,y,z\in X$ we have 
    
$$\eqalign{\phi_A(_A\langle x,y\rangle )\phi (z)&=\phi (_A\langle 
x,y\rangle\cdot z)=\phi (x\cdot\langle y,z\rangle_B)\cr
&=\phi (x)\phi_B(\langle y,z\rangle_B)=\phi (x)\langle\phi (y),\phi 
(z)\rangle_D\cr
&={}_C\langle\phi (x),\phi (y)\rangle\cdot\phi (z)\,,\cr}
$$
which implies condition (d) by Lemma \cite{fullemma}.
\eop

\uj Note that the proof of (iii) shows that if $\phi$ satisfies (b) and (c) 
then $\phi$ is an isomorphism of $X$ onto a $C-D$ Hilbert subbimodule of $
Y$.
Also note that the statements of the lemma remain true if we 
interchange condition (b) with (d) and condition (c) with (a). 

An equivalent characterization of isomorphisms between the 
imprimitivity bimodules $_AX_B$ and $_CY_D$ is a Banach space isomorphism 
$\phi :X\to Y$ satisfying the {\it ternary homomorphism identity}, that is, 
$$\phi (x\cdot\langle y,z\rangle_B)=\phi (x)\cdot\langle\phi (y),
\phi (z)\rangle_D$$
for all $x,y,z\in X$. 

\mark{idfiidlem}
\lemma$(\hbox{\rm id},\phi ,\hbox{\rm id})$ is an isomorphism between the Hilbert 
bimodules $_AX_B$ and $_AY_B$ if and only if $(\hbox{\rm id},\phi 
,\hbox{\rm id})$ is an 
isomorphism between the corresponding imprimitivity 
bimodules $_{A_0}X_{B_0}$ and $_{C_0}Y_{D_0}$.

\proof If $(\hbox{\rm id},\phi ,\hbox{\rm id})$ is an isomorphism between $_
AX_B$ and $_AY_B$ then 
$$A_0=\overline {\hbox{\rm span}}\,_A\langle X,X\rangle =\overline {\hbox{\rm span}}\,{}_
A\langle\phi (X),\phi (X)\rangle =\overline {\hbox{\rm span}}\,{}_
A\langle Y,Y\rangle =C_0\,.$$
Similarly, $B_0=D_0$ and so $(\hbox{\rm id},\phi ,\hbox{\rm id})$ is an isomorphism between $_{
A_0}X_{B_0}$ 
and $_{C_0}Y_{D_0}$.  Now suppose that $(\hbox{\rm id},\phi ,\hbox{\rm id}
)$ is an isomorphism between 
$_{A_0}X_{B_0}$ and $_{C_0}Y_{D_0}$.  If $a\in A$ and $x\in X$ then $
x=i\cdot x'$ for some $i\in A_0$ 
and $x'\in X$ and so 
$$\phi (a\cdot x)=\phi (a\cdot (i\cdot x'))=ai\cdot\phi (x')=a\cdot
\phi (i\cdot x')=a\cdot\phi (x)\,.$$
Hence $(\hbox{\rm id},\phi ,\hbox{\rm id})$ is an isomorphism between $_
AX_B$ and $_AY_B$ 
by Lemma \cite{relaxlem}.  
\eop\ 

\newsection{Morita equivalent twisted actions}

\uj Recall from [Rie] that if $_AX_B$ is an imprimitivity 
bimodule then there is a bijective correspondence (often 
called the {\it Rieffel correspondence\/}) between closed 
subbimodules of $X$ and closed ideals of $A$.  If $I$ is a 
closed ideal of $A$ then $I\cdot X$ is a closed subbimodule of $X$.  
Note that by the Cohen-Hewitt factorization theorem we 
do not have to take the closure of $I\cdot X$.  Similarly $X\cdot 
J$ 
is a closed subbimodule of $X$ if $J$ is a closed ideal of $B$.  
On the other hand if $Y$ is a closed subbimodule of $X$ then 
$_IY_J$ is an imprimitivity bimodule, where $I$ is the closed 
span of $_A$$\langle Y,Y\rangle$ and $J$ is the closed span of $\langle 
Y,Y\rangle_B$.  We 
call $_IY_J$ an {\it imprimitivity subbimodule\/} of $X$.  

\definition A {\it partial automorphism\/} of the imprimitivity 
bimodule $_AX_B$ is an isomorphism between two 
imprimitivity subbimodules of $X$.  We denote the set of 
partial automorphisms by $\hbox{\rm PAut}(X)$.  

\uj Let $A$ be a $C^{*}$-algebra, and let $S$ be a unital inverse semigroup 
with idempotent semilattice $E$, and unit $e$.  Recall from [Si2] that a 
{\it Busby-Smith twisted action\/} of $S$ on $A$ is a pair $(\beta 
,v)$, where for all 
$s\in S$, $\beta_s:A_{s^{*}}\to A_s$ is a partial automorphism, that is, an isomorphism 
between closed ideals of $A$, and for all $s,t\in S$, $v_{s,t}{}_{}$ is a unitary 
multiplier of $A_{st}$, such that for all $r,s,t\in S$ we have 

\item{(a)}{$A_e=A$;}
\item{(b)}{$\beta_s\beta_t=\hbox{\rm Ad }v_{s,t}\circ\beta_{st}$;} 
\item{(c)}{$ $$v_{s,t}=1_{M(A_{st})}$ if $s$ or $t$ is an idempotent; 
\item{(d)}{$\beta_r(av_{s,t})v_{r,st}=\beta_r(a)v_{r,s}v_{rs,t}$ for all $
a\in A_{r^{*}}A_{st}$.} 

\noindent We refer to condition (d) as the cocycle identity.

Also recall that a {\it covariant representation\/} of a Busby-Smith twisted 
action $(A,S,\beta ,v)$ is a triple $(\pi ,V,H)$, where $\pi$ 
is a nondegenerate representation of $A$ on the Hilbert 
space $H$ and $V_s$ is a partial isometry for all $s\in S$, such 
that for all $r,s\in S$ we have 

\item{(a)}{$V_s$ has initial space $\pi (A_{s^{*}})H$ and final space 
$\pi (A_s)H$;} 
\item{(b)}{$V_rV_s=\pi (v_{r,s})V_{rs}$;} 
\item{(c)}{$\pi (\beta_s(a))=V_s\pi (a)V_s^{*}\quad$for $a\in A_{
s^{*}}$.} 

\mark{morequ} 
\definition The Busby-Smith twisted 
actions $(A,S,\alpha ,u)$ and $(B,S,\beta ,v)$ are {\it Morita equivalent\/} if 
there is an imprimitivity bimodule $_A$$X_B$ and a map 
$s\mapsto (\alpha_s,\phi_s,\beta_s):S\to\hbox{\rm PAut}(X)$, such that $
\phi_s:X_{s^{*}}\to X_s$ where 
$X_s:=A_s\cdot X=X\cdot B_s$ and for all $s,t\in S$ we have 
$$\phi_s\phi_t=u_{s,t}\cdot\phi_{st}(\cdot )\cdot v_{s,t}^{*}\,.$$
We say that $(X,\phi )$ is a Morita 
equivalence between $(\alpha ,u)$ and $(\beta ,v)$, and we write 
$$(A,S,\alpha ,u)\sim_{X,\phi}(B,S,\beta ,v)\,.$$

\uj Note that $\phi_s\phi_t$ and $\phi_{st}$ have the same range $
X_{st}$ 
and so $X_{st}\subset X_s$.  

\mark{fiax}
\lemma Using the notations of Definition \cite{morequ} we 
have
\item{\rm(a)}{$\phi_s(X_{s^{*}}\cdot B_t)=X_{st};$}
\item{\rm(b)}{$\phi_s(A_{s^{*}}\cdot X_t)=X_{st};$}
\item{\rm(c)}{$\overline {\hbox{\rm span}}\,\alpha_s(_A\langle X_{
s^{*}},X_t\rangle )=A_{st},$}

\noindent for all $s,t\in S$.

\proof We know from [Si2] that $\beta_s(B_{s^{*}}B_t)=B_{st}$ and so we 
have 
$$\eqalign{\phi_s(X_{s^{*}}\cdot B_t)&=\phi_s(X_{s^{*}}\cdot B_{s^{
*}}B_t)=\phi_s(X_{s^{*}})\cdot\beta_s(B_{s^{*}}B_t)=X_s\cdot B_{s
t}=X_{st}\,,\cr}
$$
showing (a). A similar calculation shows (b). Finally (c) 
follows from the calculation:
$$\eqalign{\overline {\hbox{\rm span}}\,\alpha_s(_A\langle X_{s^{
*}},X_t\rangle )&=\overline {\hbox{\rm span}}\,\alpha_s(_A\langle 
A_{s^{*}}\cdot X_{s^{*}},X_t\rangle )=\overline {\hbox{\rm span}}
\,\alpha_s(_A\langle X_{s^{*}},A_{s^{*}}\cdot X_t\rangle )\cr
&=\overline {\hbox{\rm span}}\,{}_A\langle\phi_s(X_{s^{*}}),\phi_
s(A_{s^{*}}\cdot X_t)\rangle =\overline {\hbox{\rm span}}\,{}_A\langle 
X_s,X_{st}\rangle =A_{st}\,.\cr}
$$
\Eop\ 

\prop Morita equivalence of Busby-Smith twisted 
actions is an equivalence relation.

\proof It is easy to see that $(A,S,\alpha ,u)\sim_{A,\alpha}(A,S
,\alpha ,u)$. It 
is also easy to check that if $(A,S,\alpha ,u)\sim_{X,\phi}(B,S,\beta 
,v)$ then 
$(B,S,\beta ,v)\sim_{\tilde {X},\tilde{\phi}}(A,S,\alpha ,u)$, where $
\tilde{\phi }(\tilde {x})=\phi (x)^{\star}$. To show 
transitivity, suppose 
$$(A,S,\alpha ,u)\sim_{X,\phi}(B,S,\beta ,v)\sim_{Y,\psi}(C,S,\gamma 
,w)\,.$$
Let $Z$ be the balanced 
tensor product $X\otimes_BY$, that is, the Hausdorff completion 
of $X\odot Y$ in the $C$-valued inner product determined by 
$$\langle x_1\otimes y_1,x_2\otimes y_2\rangle_C:=\langle y_1,\langle 
x_1,x_2\rangle_B\cdot y_2\rangle_C\,.$$
It is well known that $Z$ is an $A-B$ imprimitivity bimodule.
We are going to define a map $\theta$ such that 
$(A,S,\alpha ,u)\sim_{Z,\theta}(C,S,\gamma ,w)\,.$ For all $s\in 
S$ we have 
$$\eqalign{Z_s&=(X\otimes_BY)\cdot C_s=X\otimes_B(Y\cdot C_s)\cr
&=X\otimes_B(B_s\cdot Y_s)=(X\cdot B_s)\otimes_BY_s=X_s\otimes_BY_
s\,.\cr}
$$
For all $s\in S$ the map $\theta':X_{s^{*}}\times Y_{s^{*}}\to Z_
s$ defined by 
$\theta'(x,y)=\phi_s(x)\otimes\psi_s(y)$ is bilinear, and so we have a linear map 
$\theta^{\prime\prime}_s:X_{s^{*}}\odot Y_{s^{*}}\to Z_s$ satisfying $
\theta_{s^{*}}^{\prime\prime}(x\otimes y)=\theta'(x,y)$. The 
following computation suffices to check 
that $\theta^{\prime\prime}$ is isometric:
$$\eqalign{\langle\theta_s^{\prime\prime}(x_1\otimes y_1),\theta_
s^{\prime\prime}(x_2\otimes y_2)\rangle_C&=\langle\phi_s(x_1)\otimes
\psi_s(y_1),\phi_s(x_2)\otimes\psi_s(y_2)\rangle_C\cr
&=\langle\psi_s(y_1),\langle\phi_s(x_1),\phi_s(x_2)\rangle_B\cdot
\psi_s(y_2)\rangle_C\cr
&=\langle\psi_s(y_1),\psi_s(\langle x_1,x_2\rangle_B\cdot y_2)\rangle_
C\cr
&=\gamma_s(\langle y_1,\langle x_1,x_2\rangle_B\cdot y_2\rangle_C
)\cr
&=\gamma_s(\langle x_1\otimes y_1,x_2\otimes y_2\rangle_C)\,.\cr}
$$
So $\theta_s^{\prime\prime}$ extends uniquely to an isometric linear map 
$\theta_s:Z_{s^{*}}\to Z_s$.  The above calculation also shows that $
\theta_s$ 
satisfies Definition \cite{impiso}(b), and it is routine to 
check Definition \cite{impiso}(c).  Finally for all $s,t\in S$ 
we have 
$$\eqalign{\theta_s\theta_t&=\phi_s\phi_t\otimes\psi_s\psi_t=u_{s
,t}\cdot\phi_{st}(\cdot )\cdot v_{s,t}^{*}\otimes v_{s,t}\cdot\psi_{
st}(\cdot )\cdot w_{s,t}^{*}\cr
&=u_{s,t}\cdot\phi_{st}(\cdot )\otimes v_{s,t}^{*}v_{s,t}\cdot\psi_{
st}(\cdot )\cdot w_{s,t}^{*}=u_{s,t}\cdot\theta_{st}\cdot w_{s,t}^{
*}\,.\cr}
$$
\Eop

\uj Recall [BGR] that two projections $p$ and $q$ in the 
multipliers of a $C^{*}$-algebra $C$ are called complementary 
if $p+q=1$. The two corners $pCp$ and $qCq$ are also 
called complementary.  The projection $p$ is called full if 
the corner $pCp$ is full, which means $pCp$ is not contained in 
any proper ideal of $C$ or equivalently $CpC$ is dense in $C$.  
If the $C^{*}$-algebras $A$ and $B$ are Morita equivalent then 
they are isomorphic to complementary full corners of 
the linking $C^{*}$-algebra 
$$C=\left(\matrix{A&X\cr
\tilde {X}&B\cr}
\right)$$
of $_AX_B$, where $\tilde {X}$ is the reverse module of $X$ and the 
operations on $C$ are defined by 
$$\left(\matrix{a&x\cr
\tilde {y}&b\cr}
\right)\left(\matrix{c&z\cr
\tilde {w}&d\cr}
\right)=\left(\matrix{ac+{}_A\langle x,w\rangle&a\cdot z+x\cdot d\cr
\tilde {y}\cdot c+b\cdot\tilde {w}&\langle y,z\rangle_B+bd\cr}
\right)$$
$$\left(\matrix{a&x\cr
\tilde {y}&b\cr}
\right)^{*}=\left(\matrix{a^{*}&y\cr
\tilde {x}&b^{*}\cr}
\right)\,.$$
In fact, we can identify $A$ and $B$ with $pCp$ and $qCq$ 
respectively, where 
$$p=\left(\matrix{1_{M(A)}&0\cr
0&0\cr}
\right)\quad\hbox{\rm and}\quad q=\left(\matrix{0&0\cr
0&1_{M(B)}\cr}
\right)\,.$$
Here we identified the multiplier algebra $M(C)$ with 
$$\left(\matrix{M(A)&M(X)\cr
M(\tilde {X})&M(B)\cr}
\right)$$
as in [ER, Appendix].  
On the other hand if two $C^{*}$-algebras are isomorphic to 
complementary full corners of a $C^{*}$-algebra, then they 
are Morita equivalent.

Note that if the actions $(A,S,\alpha ,u)$ and $(B,S,\beta ,w)$ are 
Morita equivalent then the $C^{*}$-algebras $A$ and $B$ are also 
Morita equivalent.  We have a natural action of $S$ on 
the linking algebra of $A$ and $B$: 

\mark{linkaction}
\prop If $(A,S,\alpha ,u)\sim_{X,\phi}(B,S,\beta ,v)$ then the formulas
$$\gamma_s\left(\matrix{a&x\cr
\tilde {y}&b\cr}
\right)=\left(\matrix{\alpha_s(a)&\phi_s(x)\cr
\phi_s(y)^{\star}&\beta_s(b)\cr}
\right),\quad w_{s,t}=\left(\matrix{u_{s,t}&0\cr
0&v_{s,t}\cr}
\right)$$
define a Busby-Smith twisted action $(C,S,\gamma ,w)$ on the 
linking algebra $C$ of $_AX_B$.  Moreover, $(Y,\gamma (\cdot )|Y)$ implements 
a Morita equivalence between $(C,S,\gamma ,w)$ and $(B,S,\beta ,v
)$, 
where $Y=\left(\smallmatrix {0&X\cr
0&B\cr}\right)\subset C\,.$

\proof It is well-known that $_CY_B$ is an imprimitivity bimodule if $
Y$ 
inherits the inner products from the $C^{*}$-algebra $C$, that 
is, $_C\langle y_1,y_2\rangle =y_1y_2^{*}$ for all $y_1,y_2\in Y$ and 
$\langle\left({0\atop 0}{x\atop b}\right),\left({0\atop 0}{z\atop 
d}\right)\rangle_B=\langle x,z\rangle_B+b^{*}d$ for all $x,z\in X$ and 
$b,d\in B$. 
It is easy to check that 
$$C_s=\left(\matrix{A_s&X_s\cr
\tilde {X}_s&B_s\cr}
\right)$$
is the closed ideal of $C$ 
which is in Rieffel correspondence with $B_s$ via the imprimitivity 
bimodule $_CY_B$.  The calculation 
$$\eqalign{\gamma_s\left(\left(\matrix{a&x\cr
\tilde {y}&b\cr}
\right)\left(\matrix{c&z\cr
\tilde {w}&d\cr}
\right)\right)&=\left(\matrix{\alpha_s(ac)+\alpha_s({}_A\langle x
,w\rangle )&\phi_s(a\cdot z+x\cdot d)\cr
\phi_s(c\cdot y+w\cdot b)^{\star}&\beta_s(\langle y,z\rangle_B+bd
)\cr}
\right)\cr
&=\left(\matrix{\alpha_s(a)&\phi_s(x)\cr
\phi_s(y)^{\star}&\beta_s(b)\cr}
\right)\left(\matrix{\alpha_s(c)&\phi_s(z)\cr
\phi_s(w)^{\star}&\beta_s(d)\cr}
\right)\cr
&=\gamma_s\left(\matrix{a&x\cr
\tilde {y}&b\cr}
\right)\gamma_s\left(\matrix{c&z\cr
\tilde {w}&d\cr}
\right)\cr}
$$
shows that $\gamma_s$ is a homomorphism for all $s\in S$.  It is 
easy to verify that $\gamma_s$ preserves adjoints and is 
bijective, hence is an isomorphism between $C_{s^{*}}$ and $C_s$.  
We only check the cocycle identity in the definition of Busby-Smith 
twisted actions.
It suffices to show that for $a\in A_{r^{*}}A_{st}$, 
$b\in B_{r^{*}}B_{st}$ and $x$, $y\in X_{r^{*}}\cap X_{st}$, 
$$\gamma_r\left(\left(\matrix{a&x\cr
\tilde {y}&b\cr}
\right)w_{s,t}\right)=\left(\matrix{\alpha_r(au_{s,t})u_{r,st}&\phi_
r(x\cdot v_{s,t})v_{r,st}\cr
\phi_r(u_{s,t}\cdot y)^{\star}\cdot u_{r,st}&\beta_r(bv_{s,t})v_{
r,st}\cr}
\right)$$
and
$$\gamma_r\left(\matrix{a&x\cr
y^*&b\cr}
\right)w_{r,s}w_{rs,t}=\left(\matrix{\alpha_r(a)u_{r,s}u_{rs,t}&\phi_
r(x)v_{r,s}v_{rs,t}\cr
\phi_r(y)^\star\cdot u_{r,s}u_{rs,t}&\beta_r(b)v_{r,s}v_{rs,t}\cr}
\right)$$
are the same. The diagonals are clearly equal. We check 
the upper right corner. Since $x=x_r\cdot a_r$ for some $x_r\in X_{
r^{*}}$ and 
$a_r\in A_{r^{*}}$ we have
$$\eqalign{\phi_r(xv_{s,t})v_{r,st}&=\phi_r(x_r\cdot a_rv_{s,t})=
\phi_r(x_r)\beta_r(a_rv_{s,t})\cr
&=\phi_r(x_r)\beta_r(a_r)v_{r,s}v_{rs,t}=\phi_r(x)v_{r,s}v_{rs,t}\,
.\cr}
$$
The equality of the lower left corners follows similarly.  For the 
other part, the 
conditions of Definition \cite{morequ} for the pair $(Y,\gamma (\cdot 
)|Y)$ follow 
from routine calculations. 
\eop

\uj A similar proof shows that in the previous theorem 
$(C,S,\gamma ,w)$ and $(A,S,\beta ,v)$ are also Morita equivalent.
Recall [Si2] that two Busby-Smith twisted actions $(\alpha ,u)$ 
and $(\beta ,w)$ of $S$ on $A$ are {\it exterior equivalent\/} if for all 
$s\in S$ there is a unitary multiplier $V_s$ of $E_s$ such that 
for all $s,t\in S$ 
\item{(a)}{$\beta_s=\hbox{\rm Ad }V_s\circ\alpha_s$;} 
\item{(b)}{$w_{s,t}=V_s\alpha_s(1_{M(E_{s^{*}})}V_t)u_{s,t}V_{st}^{
*}$.} 

\theorem If the twisted actions $(A,S,\alpha ,u)$ and $(A,S,\beta 
,w)$ are 
exterior equivalent, then they are also Morita 
equivalent.

\proof Let $V$ implement an exterior equivalence between 
$(\alpha ,u)$ and $(\beta ,w)$.  We show that $(A,\phi )$ implements the 
Morita equivalence, where $\phi_s:A_{s^{*}}\to A_s$ is defined by 
$\phi_s(a)=\alpha_s(a)V_s^{*}$. For $a,b,x\in A_{s^{*}}$ we have 
$$\eqalign{\phi_s(x\cdot b)&=\alpha_s(x)\alpha_s(b)V_s^{*}=\alpha_
s(x)V_s^{*}\beta_s(b)=\phi_s(x)\cdot\phi_s(x)\cr}
$$
verifying Definition \cite{impiso}(a). 
If $x,y\in X_{s^{*}}=A_{s^{*}}$, then we have 
$$\alpha_s(_A\langle x,y\rangle )=\alpha_s(xy^{*})=\alpha_s(x)V^{
*}(\alpha_s(y)V^{*})^{*}={}_A\langle\phi_s(x),\phi_s(y)\rangle ,$$
which verifies Condition \cite{impiso}(d).  By the note after Lemma 
\cite{relaxlem}, it remains to observe that if $x\in X_{(st)^{*}}
=A_{(st)^{*}}$ then 
$$\eqalign{\big(\phi_s\phi_t\big)(x)&=\alpha_s(\alpha_t(x)V_t^{*}
)V_s^{*}\cr
&=\alpha_s(\alpha_t(x))u_{s,t}V_{st}^{*}V_{st}u_{s,t}^{*}\alpha_s
(1_{M(A_{s^{*}})}V_t)^{*}V_s^{*}\cr
&=u_{s,t}\phi_{st}(x)w_{s,t}^{*}\,.\cr}
$$
\Eop

\uj Recall [Rie] that if $_AX_B$ is an imprimitivity bimodule then every 
representation $\pi$ of $B$ on a Hilbert space $H$ induces a representation 
$\pi^X$ of $A$ on the Hilbert space $H^X$ defined by $\pi^X(a)(x\otimes
\xi )=(a\cdot x)\otimes\xi$, 
where $H^X$ is the Hausdorff completion $X\otimes_BH$ of the algebraic tensor 
product $X\odot H$ in the seminorm generated by the semi-inner product 
$$(x\otimes\xi\mid y\otimes\eta ):=(\pi (\langle y,x\rangle_B)\xi
\mid\eta )_H=(\xi\mid\pi (\langle x,y\rangle_B)\eta ))_H\,.$$
Note that $(x\cdot b)\otimes\xi =x\otimes\pi (b)\xi$ for all $x\in 
X$, $b\in B$ and $\xi\in H$.  The 
following is the semigroup version of [Com, Section 2].  

\theorem If $(A,S,\alpha ,u)\sim_{X,\phi}(B,S,\beta ,v)$ then every covariant 
representation $(\pi ,V,H)$ of $(\beta ,v)$ induces a covariant 
representation $(\pi^X,V^X,H^X)$ of $(\alpha ,u)$ on the Hilbert space 
$H^X=X\otimes_BH$, where $\pi^X$ is as above and
$$V_s^X(x\otimes\xi )=\phi_s(x)\otimes V_s(\xi )$$
for all elementary tensors $x\otimes\xi\in H_{s^{*}}^X=X_{s^{*}}\otimes_
BH$.

\proof First note that if $x\in X_s$ and $\xi\in H$ then $x=y\cdot 
b$ for some 
$y\in X_s$ and $b\in B_s$, hence $x\otimes\xi =(y\cdot b)\otimes\xi 
=y\otimes\pi (b)\xi$.  So 
$H_s^X=X_s\otimes_BH_s$ where $H_s=\pi (B_s)H=V_sH$.  To show 
the existence of $V_s^X$, define $T:X_{s^{*}}\times H\to X_s\otimes_
BH_s$ by 
$T(x,\xi )=\phi_s(x)\otimes V_s\xi$.  $T$ is clearly bilinear so there is a 
unique linear map $T':X_{s^{*}}\odot H_{s^{*}}\to X_{s^{*}}\otimes_
BH_{s^{*}}$ such that 
$T'(x\otimes\xi )=T(x,\xi )$.  We check that $T'$ is isometric.  For 
$x,y\in X_{s^{*}}$ and $\xi ,\eta\in H_s$ we have 
$$\eqalign{(T'(x\otimes\xi )\mid T'(y\otimes\eta ))_{H^X}&=(\phi_
s(x)\otimes V_s\xi\mid\phi_s(y)\otimes V_s\eta )_{H^X}\cr
&=(\pi (\langle\phi_s(y),\phi_s(x)\rangle_B)V_s\xi\mid V_s\eta )_
H\cr
&=(\pi (\beta_s(\langle y,x\rangle_B))V_s\xi\mid V_s\eta )_H\cr
&=(V_s\pi (\langle y,x\rangle_B)\xi\mid V_s\eta )_H\cr
&=(\pi (\langle y,x\rangle_B)\xi\mid\eta )_H=(x\otimes\xi\mid y\otimes
\eta )_{H^X}\,.\cr}
$$
So $T'$ determines an isometry $T^{\prime\prime}$ from $H_{s^{*}}^
X$ to $H_s^X$.  If we define 
$V_s^X$ to be $T^{\prime\prime}$ on $H_{s^{*}}^X$ and $0$ on $(H_{
s^{*}}^X)^{\perp}$ then $V_s^X$ is a partial 
isometry with initial space 
$H_{s^{*}}^X=(A_{s^{*}}\cdot X)\otimes_BH=\pi^X(A_{s^{*}})H^X$.  It follows that the 
final space of $V_s^X$ is $\pi^X(A_s)H^X$.  

We can check the covariance condition for elementary tensors.  Let 
$a\in A_{s^{*}}$ and $x\otimes\xi\in X\odot H$.  Since $H=H_s\oplus 
H_s^{\perp}$, we only need to consider 
the two cases $\xi\in H_s$ and $\xi\in H_s^{\perp}$.  If $\xi\in 
H_s$ then $\xi =\pi (ab)\eta$ for some 
$a,b\in A_s$ and $\eta\in H$.  Hence $x\otimes\xi =x\cdot a\otimes
\pi (b)\eta$ and so we can assume 
that $x\in X_s$.  Thus, 
$$\eqalign{V_s^X\pi^X(a)(V_s^X)^{*}(x\otimes\xi )&=\phi_s(a\cdot\phi_
s^{-1}(x))\otimes V_sV_s^{*}(\xi ))\cr
&=(\alpha_s(a)\cdot x)\otimes\xi =\pi^X(\alpha_s(a))(x\otimes\xi 
)\,.\cr}
$$
On the other hand if $\xi\in (H_s)^{\perp}$ then for all $y\in X_
s$ and $\eta\in H$ we have
$$\eqalign{(x\otimes\xi\mid y\otimes_B\eta )_H=(\xi\mid\pi (\langle 
x,y\rangle_B)\eta )_H=0\cr}
$$
and so $x\otimes\xi\in (H_s^X)^{\perp}$. This means $(V_s^X)^{*}(
x\otimes\xi )=0$. 
Since $\alpha_s(a)\cdot x$ is in $X_s$ it is of the form $y\cdot 
b$ for some $y\in X$ and 
$b\in B_s$. Thus,
$$\pi^X(\alpha_s(a))(x\otimes\xi )=(\alpha_s(a)\cdot x)\otimes\xi 
=y\otimes\pi (b)\xi =0$$
as well.
\Eop

\uj Of course the inducing process works the other way 
too, that is, every covariant representation of $\alpha$ induces 
a covariant representation of $\beta$. 

Recall [Si2] that the crossed product $A\times_{\alpha ,u}S$ of a 
Busby-Smith twisted action $(A,S,\alpha ,u)$ is the Hausdorff completion 
of the Banach $*$-algebra 
$$L_{\alpha}=\{x\in l^1(S,A):x(s)\in A_s\hbox{\rm \ for all }s\in 
S\}$$
with operations
$$\eqalign{\big(x*y\big)(s)&=\sum_{rt=s}\alpha_r(\alpha_r^{-1}(x(
r))y(t))u_{r,t}\quad\hbox{\rm and}\quad x^{*}(s)=u_{s,s^{*}}\alpha_
s(x(s^{*})^{*})\cr}
$$
in the $C^{*}$-seminorm $\|\cdot\|_{\alpha}$ defined by 
$$\|x\|_{\alpha}=\sup\{\|(\pi\times V)(x)\|:(\pi ,V)\hbox{\rm \ is a covariant representation of }
(A,S,\alpha ,u)\}.$$
Alternatively, generalizing Paterson's approach [Pat] to the twisted 
case, we could define  
$$\|x\|_{\alpha}=\sup\{\|\phi (x)\|:\phi\hbox{\rm \ is a coherent representation of }
L_{\alpha}\}$$
where a representation $\phi$ of $L_{\alpha}$ is coherent if it satisfies 
$\phi (a\chi_{ss^{*}})=\phi (a\chi_e)$ for all $s\in S$.  So if $
I_{\alpha}$ is the closed ideal 
generated by elements of the form $a\chi_{ss^{*}}-a\chi_e$, then the crossed 
product $A\times_{\alpha}S$ is the enveloping $C^{*}$-algebra of $
L_{\alpha}/I_{\alpha}$.  

If $\chi_s$ denotes the characteristic function of $\{s\}$, then $
a\chi_s$ 
is an element of $L_{\alpha}$ for all $a\in A_s$. The canonical image 
of $a\chi_s$ in $A\times_{\alpha ,u}S$ will be denoted by $a\delta_
s$. Then 
$A\times_{\alpha ,u}S$ is the closed span of $\{a\delta_s:a\in A_
s,\,s\in S\}$. Note 
that we have the following formulas:
$$\eqalign{a_s\delta_s*a_t\delta_t&=\alpha_s(\alpha_s^{-1}(a_s)a_
t)u_{s,t}\delta_{st}\cr
(a\delta_s)^{*}&=\alpha_s^{-1}(a^{*})u_{s^{*},s}^{*}\delta_{s^{*}}
.\cr}
$$

\uj The idea of the proof of the following theorem comes 
from [Com], [CMW] and [Kal].

\mark{moreqbsxp} \theorem If $(A,S,\alpha ,u)$ and $(B,S,\beta ,v
)$ are 
Morita equivalent actions, then the crossed products 
$A\times_{\alpha ,u}S$ and $B\times_{\beta ,v}S$ are also Morita equivalent.  

\proof Let $(X,\phi )$ be a Morita equivalence, and let $(\gamma 
,w)$ be the 
Busby-Smith twisted action of $S$ on the linking algebra $C$ of $_
AX_B$ as 
in Proposition \cite{linkaction}.  It suffices to show that $A\times_{
\alpha ,u}S$ and 
$B\times_{\beta ,v}S$ are complementary full corners of $C\times_{
\gamma ,w}S$.  Let 
$$p=\left(\matrix{1_{M(A)}&0\cr
0&0\cr}
\right)\quad\hbox{\rm and}\quad q=\left(\matrix{0&0\cr
0&1_{M(B)}\cr}
\right)\,.$$
It is clear that $p\delta_e$ and $q\delta_e$ are complementary projections in 
$M(C\times_{\gamma ,w}S)$.  We show that $q\delta_e$ is a full projection.  If 
$$c=\left(\matrix{a_s&x_s\cr
\tilde {y}_s&b_s\cr}
\right)\in C_s\quad\hbox{\rm and}\quad d=\left(\matrix{a_t&x_t\cr
\tilde {y}_t&b_t\cr}
\right)\in C_t$$
then 
$$c\delta_s*q\delta_e*d\delta_t=\left(\matrix{u_{s,t}\alpha_s(_A\langle
\phi_s^{-1}(x_s),y_t\rangle )&\phi_s(\phi_s^{-1}(x_s)\cdot b_t)\cdot 
v_{s,t}\cr
\phi_s(\beta_s^{-1}(b_s)\cdot y_t)^{\star}\cdot u_{s,t}&\beta_s(\beta_
s^{-1}(b_s)b_t)v_{s,t}\cr}
\right)\delta_{st}\,.$$
We can check fullness on the four corners, and this 
can be done easily using Lemma \cite{fiax}. A similar 
calculation shows that $p\delta_e$ is also full.  

Now we show that $B\times_{\beta ,v}S=q\delta_e*(C\times_{\gamma 
,w}S)*q\delta_e$.  We use the fact 
that $B\times_{\beta ,v}S$ is the Hausdorff completion $\overline {
L_{\beta}}^{\,\|\cdot\|_{\beta}}$ of $L_{\beta}$ in the greatest 
$C^{*}$-seminorm $\|\cdot\|_{\beta}$ coming from covariant representations of $
(\beta ,v)$, 
while $C\times_{\gamma ,w}S$ is the Hausdorff completion $\overline {
L_{\gamma}}^{\,\|\cdot\|_{\gamma}}$ of $L_{\gamma}$ in the 
greatest $C^{*}$-seminorm $\|\cdot\|_{\gamma}$ coming from covariant representations of 
$(\gamma ,w)$.  Since 
$$q\delta_e(C\times_{\gamma ,w}S)q\delta_e=q\delta_e(\overline {L_{
\gamma}}^{\,\|\cdot\|_{\gamma}})q\delta_e=\overline {q\chi_e*L_{\gamma}
*q\chi_e}^{\,\|\cdot\|_{\gamma}}=\overline {L_{\beta}}^{\,\|\cdot
\|_{\gamma}}\,,$$
it suffices to show that the seminorms $\|\cdot\|_{\beta}$ and $\|
\cdot\|_{\gamma}$ 
are the same on $L_{\beta}$, where we regard $L_{\beta}$ as a 
subspace of $L_{\gamma}$.  If $(\pi ,V)$ is a covariant representation of 
$(\gamma ,w)$ then $(\pi |B,\pi (q)V)$ is a covariant representation of 
$(\beta ,v)$ and so $\|\cdot\|_{\gamma}\leq\|\cdot\|_{\beta}$ on $
L_{\beta}$.  On the other hand, a 
covariant representation $(\pi ,V,H)$ of $(\beta ,v)$ induces a 
covariant representation $(\pi^Y,V^Y,H^Y)$ of $(\gamma ,w)$, where 
$Y=\left(\smallmatrix {0&X\cr
0&B\cr
}\right)$, $H^Y=Y$$\otimes_BH$ and 
$$\pi^Y\left(\matrix{a&x\cr
\tilde {y}&b\cr}
\right)\left(\left(\matrix{0&z\cr
0&d\cr}
\right)\otimes\xi\right)=\left(\matrix{0&a\cdot z+x\cdot d\cr
0&\langle y,z\rangle_B+bd\cr}
\right)\otimes\xi\,,$$
$$V_s^Y\left(\left(\matrix{0&z\cr
0&d\cr}
\right)\otimes\xi\right)=\left(\matrix{0&\phi_s(z)\cr
0&\beta_s(d)\cr}
\right)\otimes V_s\xi\,.$$
The image of $\left(\smallmatrix{0&0\cr
0&b\cr}
\right)\chi_s\in L_{\beta}$ under $\pi^Y\times V^Y$ evaluated at an 
elementary tensor $\left(\smallmatrix {0&z\cr
0&d\cr}\right)\otimes\xi$ of $H^Y$ is 
$$\eqalign{\pi^Y\left(\matrix{0&0\cr
0&b\cr}
\right)V_s^Y\left(\left(\matrix{0&z\cr
0&d\cr}
\right)\otimes\xi\right)&=\left(\matrix{0&0\cr
0&b\beta_s(d)\cr}
\right)\otimes V_s\xi\cr
&=\left(\smallmatrix{0&0\cr
0&b\cr}
\right)\otimes\pi (\beta_s(d))V_s\xi\cr
&=\left(\smallmatrix{0&0\cr
0&b\cr}
\right)\otimes V_s\pi (d)\xi\,.\cr}
$$
If $d\in B$ and $\xi\in H$ then
$$\eqalign{\|\left(\smallmatrix {0&0\cr
0&d\cr}\right)\otimes\xi\|_{H^Y}^2&=\big(\left(\smallmatrix {0&0\cr
0&d\cr}\right)\otimes\xi\mid\left(\smallmatrix {0&0\cr
0&d\cr}\right)\otimes\xi\big)_{H^Y}\cr
&=(\pi (\langle\left(\smallmatrix {0&0\cr
0&d\cr}\right),\left(\smallmatrix {0&0\cr
0&d\cr}\right)\rangle_B)\xi\mid\xi )_H\cr
&=(\pi (d^{*}d)\xi\mid\xi )_H=(\pi (d)\xi\mid\pi (d)\xi )_H\cr
&=\|\pi (d)\xi\|_H^2\,.\cr}
$$
Hence if $b_i\in B_{s_i}$ for all $i=1,\dots,n$ and 
$f=\sum_{i=1}^nb_i\chi_{s_i}\in L_{\beta}$ then  
$$\|\pi^Y\times V^Y(f)\big(\left(\smallmatrix {0&0\cr
0&d\cr}\right)\otimes\xi\big)\|_{H^Y}=\|\pi\times V(f)\pi (d)\xi\|_
H\,.$$
On the other hand
$$\eqalign{\|\pi\times V(f)\|&=\sup\{\,{{\|\pi\times V(f)\pi (d)\xi
\|_H}\over {\|\pi (d)\xi\|_H}}:d\in B,\xi\in H\}\cr
&=\sup\{\,{{\|\pi^Y\times V^Y(f)\big(\left(\smallmatrix {0&0\cr
0&d\cr
}\right)\otimes\xi\big)\|_{H^Y}}\over {\|\left(\smallmatrix {0&0\cr
0&d\cr
}\right)\otimes\xi\|_{H^Y}}}:d\in B,\xi\in H\}\cr
&\leq\|\pi^Y\times V^Y(f)\|\cr}
$$
which implies that $\|\cdot\|_{\gamma}\ge\|\cdot\|_{\beta}$ on $L_{
\beta}$.
A similar argument shows that $A\times_{\alpha ,u}S=p\delta_e*(C\times_{
\gamma ,w}S)*p\delta_e$.
\eop

\uj The proof also shows that if we use the notation 
$X\times_{u,\phi ,v}S:=p\delta_e(C\times_{\gamma ,w}S)q\delta_e$ or simply $
X\times S$, then 
$$A\times_{\alpha ,u}S\sim_{X\times S}B\times_{\beta ,v}S\,.$$

We now have two different ways to induce 
representations of $A\times_{\alpha}S$ from representations of $B
\times_{\beta}S$.  
The next result shows that they are essentially the 
same. For simplicity we only state the untwisted version of the 
result because that is all we need later.  
The proof closely follows that of similar results 
in [Ech] and [Kal], and goes back ultimately to [Com].  

\mark{unequ}
\prop If $(A,S,\alpha )\sim_{X,\phi}(B,S,\beta )$ and $\pi\times 
V$ is a representation 
of $B\times_{\beta}S$ on $H$, then the induced representations $\pi^
X\times V^X$ and 
$(\pi\times V)^{X\times_{\phi}S}$ are unitarily equivalent.

\proof Let $Y=X\times_{\phi}S$. The map 
$$T':X\times H\to H^Y\quad\hbox{\rm defined by}\quad T'(x,\xi )=
\left(\smallmatrix{0&x\cr
0&0\cr}
\right)\delta_e\otimes\xi$$
is bilinear and so there is a unique linear map $T^{\prime\prime}
:X\odot H\to H^Y$ 
such that $T'{}'(x\otimes\xi )=T'(x,\xi )$. We check that $T^{\prime
\prime}$ is 
isometric. For $x,y\in X$ and $\xi ,\eta\in H$ we have  
$$\eqalign{(T^{\prime\prime}(x\otimes\xi )\mid T^{\prime\prime}(y
\otimes\eta ))_{H^Y}&=\left(\left(\smallmatrix{0&x\cr
0&0\cr}
\right)\delta_e\otimes\xi\mid\left(\smallmatrix{0&y\cr
0&0\cr}
\right)\delta_e\otimes\eta\right)_{H^Y}\cr
&=\left((\pi\times V)\left(\left<\left(\smallmatrix{0&y\cr
0&0\cr}
\right)\delta_e,\left(\smallmatrix{0&x\cr
0&0\cr}
\right)\delta_e\right>_{B\times_{\beta}S}\right)\xi\mid\eta\right)_H\cr
&=\big(\big(\pi\times V\big)(\langle y,x\rangle_B\delta_e)\xi\mid\eta 
\big)_H=\big(\pi (\langle y,x,\rangle_B)\xi\mid\eta \big)_H\cr
&=(x\otimes\xi\mid y\otimes\eta )_{H^X}\,.\cr}
$$
So we have an isometry $T:H^X\to H^Y$ such that 
$T(x\otimes\xi )=\left(\smallmatrix{0&x\cr
0&0\cr}
\right)\delta_e\otimes\xi$. $T$ is onto since if $x_s\in X_s$ and $
\xi\in H$ 
then there are $x\in X$ and $b\in B_s$ such that $x_s=x\cdot b$ and 
so
$$\eqalign{\left(\matrix{0&x_s\cr
0&0\cr}
\right)\delta_s\otimes\xi&=\left(\matrix{0&x\cr
0&0\cr}
\right)\delta_e\left(\matrix{0&0\cr
0&b\cr}
\right)\delta_s\otimes\xi =\left(\left(\matrix{0&x\cr
0&0\cr}
\right)\delta_e\cdot b\delta_s\right)\otimes\xi\cr
&=\left(\matrix{0&x\cr
0&0\cr}
\right)\delta_e\otimes\pi\times V(b\delta_s)\xi =T(x\otimes\pi (b
)V_s\xi )\,.\cr}
$$
Carrying the above calculation a little further, we have 
$$\left(\matrix{0&x_s\cr
0&0\cr}
\right)\delta_s\otimes\xi =T((x\cdot b)\otimes V_s\xi )=\left(\matrix{
0&x_s\cr
0&0\cr}
\right)\delta_e\otimes V_s\xi\,,$$
and we need this fact in the verification that $T$ intertwines $\pi^
X\times V^X$ 
and $(\pi\times V)^Y$:  for $a\in A_s$ and $x\otimes\xi\in X\odot 
H$ we have 
$$\eqalign{&(\pi\times V)^Y\big(a\delta_s\big)T(x\otimes\xi )=\Big
(a\delta_s\cdot\left(\smallmatrix {0&x\cr
0&0\cr}\right)\delta_e\Big)\otimes\xi\cr
&\quad =\left(\left(\smallmatrix {a&0\cr
0&0\cr}\right)\delta_s\left(\smallmatrix {0&x\cr
0&0\cr}\right)\delta_e\right)\otimes\xi =\gamma_s\left(\gamma_s^{
-1}\left(\smallmatrix {
a&0\cr
0&0\cr}\right)\left(\smallmatrix {0&x\cr
0&0\cr}\right)\right)\delta_s\otimes\xi\cr
&\quad =\left(\matrix{0&\phi_s(\alpha_s^{-1}(a)\cdot x)\cr
0&0\cr}
\right)\delta_s\otimes\xi =\left(\matrix{0&\phi_s(\alpha_s^{-1}(a
)\cdot x)\cr
0&0\cr}
\right)\delta_e\otimes V_s\xi\cr
&\quad =T(\phi_s(\alpha_s^{-1}(a)\cdot x)\otimes V_s\xi )=TV_s^X\pi^
X(\alpha_s^{-1}(a))(x\otimes\xi )\cr
&\quad =T\pi^X(a)V_s^X(x\otimes\xi )=T(\pi^X\times V^X)\big(a\delta_
s\big)(x\otimes\xi )\,.\cr}
$$
\Eop

\uj Let $A$ be a $C^{*}$-algebra, let $S$ be a unital inverse semigroup with 
idempotent semilattice $E$, and let $N$ be a 
{\it normal Clifford subsemigroup\/} of $S$.  Recall from [Si2] that a 
subsemigroup $N$ of $S$ is a normal Clifford subsemigroup if it is 
normal, that is, $E\subset N$ and $sNs^{*}$$\subset N$, and it is also Clifford, that is, 
$n^{*}n=nn^{*}$ for all $n\in N$.  
Also recall from [Si2] that a {\it Green twisted action\/} of $(S
,N)$ on $A$ is a 
pair $(\gamma ,\tau )$, where $\gamma$ is an inverse semigroup action of $
S$ on $A$ (that 
is, a semigroup homomorphism $s\mapsto (\gamma_s,A_{s^{*}},A_s):S
\to\hbox{\rm PAut}\,(A)$ with 
$A_e=A$) and $\tau_n$ is a unitary multiplier of $A_n$ for all $n
\in N$, such that 
for all $n,l\in N$ we have 
\item{(a)}{$\gamma_n=\hbox{\rm Ad }\tau_n$;} 
\item{(b)}{$\gamma_s(\tau_n)=\tau_{sns^{*}}$ for all $s\in S$ with $
n^{*}n\leq s^{*}s$;} 
\item{(c)}{$\tau_n\tau_l=\tau_{nl}$.} 

\noindent The following is the semigroup version of [Ech, Definition 1].  

\definition The Green twisted actions $(A,S,N,\alpha ,\tau )$ and $
(B,S,N,\beta ,\rho )$ are 
{\it Morita equivalent\/} if there is a Morita equivalence $(X,\phi 
)$ between the 
untwisted actions $(A,S,\alpha )$ and $(B,S,\beta )$ such that $\tau_
n\cdot x=\phi_n(x)\cdot\rho_n$ for 
all $n\in N$ and $x\in X_n$.  We say that $(X,\phi )$ is a {\it Morita equivalence \/}
between $(A,S,N,\alpha ,\tau )$ and $(B,S,N,\beta ,\rho )$, and we write 
$(A,S,N,\alpha ,\tau )\sim_{X,\phi}(A,S,N,\beta ,\rho )$.  

\uj The proof of the following theorem is modeled on Echterhoff's 
proof [Ech] in the group case.  

\theorem If $(A,S,N,\alpha ,\tau )$ and $(B,S,N,\beta ,\rho )$ are Morita 
equivalent Green twisted actions then the crossed products 
$A\times_{\alpha ,\tau}S$ and $A\times_{\beta ,\rho}S$ are also Morita equivalent.

\proof Let $(X,\phi )$ be a Morita equivalence.  Suppose $(\pi ,V
,H)$ is 
a covariant representation of $\beta$ which preserves the twist, that is, 
$\pi (\rho_n)=V_n$ for all $n\in N$.  The induced representation $
(\pi^X,V^X,H^X)$ of $\alpha$ 
also preserves the twist, since if $x,y\in X_n$ and $\xi ,\eta\in 
H_n$ then 
$$\eqalign{(\pi^X(\tau_n)(x\otimes\xi )\mid y\otimes\eta )_{H^X}&
=(\tau_n\cdot x\otimes\xi\mid y\otimes\eta )_{H^X}=(\pi (\langle 
y,\tau_n\cdot x\rangle_B)\xi\mid\eta )_H\cr
&=(\pi (\langle y,\phi_n(x)\rangle_B\rho_n)\xi\mid\eta )_H=(\pi (
\langle y,\phi_n(x)\rangle_B)V_n\xi\mid\eta )_H\cr
&=(V_n^X(x\otimes\xi )\mid y\otimes\eta )_{H^X}\cr}
$$
and so $\pi^X(\tau_n)=V_n^X$. A similar calculation shows that if 
$(\pi^X,V^X)$ preserves the twist then so does $(\pi ,V)$.  By [Rie, 
Proposition 3.3] the kernels of $\pi\times V$ and $(\pi\times V)^{
X\times_{\phi}S}$ are 
in Rieffel correspondence.  By Proposition \cite{unequ}, 
$\pi^X\times V^X$ and $(\pi\times V)^{X\times_{\phi}S}$ have the same kernel.  Hence the 
twisting ideals of $I_{\tau}$ and $I_{\rho}$ are in Rieffel 
correspondence and so the quotients are Morita 
equivalent by [Rie, Corollary 3.2].  

\newsection{Connection with twisted partial actions}

\uj The close connection between partial actions and 
inverse semigroup actions [Si1], [Ex3], [Si2] makes it 
possible to get quick results about the Morita 
equivalence of crossed products of twisted partial 
actions. First recall the definition of a twisted partial action from 
[Ex2].

\mark{partialaction}
\definition A (discrete) {\it twisted partial action\/} of a group $
G$ on a 
$C^{*}$-algebra $A$ is a pair $(\alpha ,u)$, where for all $s\in 
G$, 
$\alpha_s:A_{s^{-1}}\to A_s$ is a partial automorphism of $A$, and for all 
$r$, $s\in G$, $u_{r,s}$ is a unitary multiplier of $A_rA_{rs}$, such that 
for all $r$, $s$, $t\in G$ we have 
\item{(a)}{$A_e=A$, and $\alpha_e$ is the identity automorphism 
of $A$;} 
\item{(b)}{$\alpha_r(A_{r^{-1}}A_s)=A_rA_{rs}$;} 
\item{(c)}{$\alpha_r(\alpha_s(a))=u_{r,s}\alpha_{rs}(a)u_{r,s}^{*}$ for all 
$a\in A_{s^{-1}}A_{s^{-1}r^{-1}}$;} 
\item{(d)}{$u_{e,t}=u_{t,e}=1_{M(A)}$;} 
\item{(e)}{$\alpha_r(au_{s,t})u_{r,st}=\alpha_r(a)u_{r,s}u_{rs,t}$ for all 
$a\in A_{r^{-1}}A_sA_{st}$;}  

\mark{moreqpar} \definition The twisted partial actions 
$(A,G,\alpha ,u)$ and $(B,G,\mu ,w)$ are {\it Morita equivalent\/} if there 
is an imprimitivity bimodule $_AX_B$ and a map 
$s\mapsto (\alpha_s,\phi_s,\mu_s):G\to\hbox{\rm PAut}\,(X)$, such that $
\phi_s:X_{s^{*}}\to X_s$ where 
$X_s=A_s\cdot X=X\cdot B_s$ and for all $s,t\in G$ and 
$x\in X\cdot B_{t^{-1}s^{-1}}B_{t^{-1}}$ we have 
$$\phi_s\phi_t(x)=u_{s,t}\cdot\phi_{st}(x)\cdot w_{s,t}^{*}\,.$$
We say that $(X,\phi )$ is a Morita equivalence between $
(A,G,\alpha ,u)$ and 
$(B,G,\theta ,w)$, and we write 
$$(A,G,\alpha ,u)\sim_{X,\phi}(B,G,\mu ,w)\,.$$

\uj Recall from [Ex3] that for a group $G$, the {\it associated }
{\it inverse semigroup} $S(G)$ has elements written in  
canonical form $[g_1][g_1^{-1}]\cdots [g_m][g_m^{-1}][s]$, where 
$g_1,\dots,g_n,s\in G$, and the order of the $[g_i][g_i^{-1}]$ terms is 
irrelevant. Multiplication and inverses are defined by
$$\eqalign{&[g_1][g_1^{-1}]\cdots [g_m][g_m^{-1}][s]\cdot [h_1][h_
1^{-1}]\cdots [h_m][h_m^{-1}][t]\cr
&\quad =[g_1][g_1^{-1}]\cdots [g_m][g_m^{-1}][s][s^{-1}][sh_1][(s
h_1)^{-1}]\cdots [sh_m][(sh_m)^{-1}][st]\cr}
$$
and
$$([g_1][g_1^{-1}]\cdots [g_m][g_m^{-1}][s])^{*}=[s^{-1}g_m][(s^{
-1}g_m)^{-1}]\cdots [s^{-1}g_1][(s^{-1}g_1)^{-1}][s^{-1}]\,.$$
Thus [e] is an identity element for $S(G)$ if $e$ is the identity of $
G$, so 
we can write $[g_1][g_1^{-1}]\cdots [g_m][g_m^{-1}]$ for $[g_1][g_
1^{-1}]\cdots [g_m][g_m^{-1}][e]$, and these 
are the idempotents of $S(G)$.
Recall from [Si2, Section 4] that if $(A,G,\alpha ,u)$ is a twisted partial 
action, then the corresponding Busby-Smith twisted action 
$(A,S(G),\beta ,v)$ is defined by
 
$$A_p=A_{g_1}\cdots A_{g_m}A_s$$
$$\beta_p=\alpha_{g_1}\alpha_{g_1}^{-1}\cdots\alpha_{g_m}\alpha_{
g_m}^{-1}\alpha_s\,$$
$$v_{p,q}=1_{M(A_{pq})}u_{s,t}\,,$$
where 
$$p=[g_1][g_1^{-1}]\cdots [g_m][g_m^{-1}][s],\quad q=[h_1][h_1^{-
1}]\cdots [h_n][h_n^{-1}][t]\,.$$

\mark{moreqcor}
\theorem The twisted partial actions $(A,G,\alpha ,u)$ and 
$(B,G,\mu ,w)$ are Morita equivalent if and only if the 
corresponding Busby-Smith twisted actions $(A,S(G),\beta ,v)$ 
and $(B,S(G),\nu ,z)$ are Morita equivalent. 

\proof Suppose $(A,S(G),\beta ,v)\sim_{X,\phi}(A,S(G),\nu ,z)$.  If we 
identify the element $s\in G$ with $[s]\in S(G)$, then 
$\phi :G\to\hbox{\rm PAut}\,(X)$. For $s,t\in G$ and 
$x\in X\cdot B_{t^{-1}s^{-1}}B_{t^{-1}}$ we have 
$$\eqalign{\phi_s\phi_t(x)&=\phi_{[s]}\phi_{[t]}(x)=v_{[s],[t]}\cdot
\phi_{[s][t]}(x)\cdot z_{[s],[t]}^{*}\cr
&=v_{[s],[t]}\cdot\phi_{[st][t^{-1}][t]}(x)\cdot z_{[s],[t]}^{*}\cr
&=v_{[s],[t]}v_{[st],[t^{-1}][t]}^{*}\cdot\phi_{[st]}\phi_{[t^{-1}
][t]}(x)\cdot z_{[st],[t^{-1}][t]}z_{[s],[t]}^{*}\cr
&=u_{s,t}\cdot\phi_{[st]}(x)\cdot w_{s,t}^{*}\qquad\hbox{\rm since, e.g., $
v_{[st],[t^{-1}][t]}=1_{M(A_{[s][t]})}$}\cr
&=u_{s,t}\cdot\phi_{st}(x)\cdot w_{s,t}^{*}\,.\cr}
$$
Now suppose $(A,G,\alpha ,u)\sim_{X,\phi}(B,G,\mu ,w)$.  We can extend $
\phi$ 
to $S(G)$ by defining 
$$\phi_p=\phi_{g_1}\phi_{g_1}^{-1}\cdots\phi_{g_m}\phi_{g_m^{-1}}
\phi_s$$
for $p=[g_1][g_1^{-1}]\cdots [g_m][g_m^{-1}][s]\in S(G)$. We verify 
Definition \cite{impiso}(d). If 
$$x,y\in X_{p^{*}}=X\cdot B_{p^{*}}=B_sB_{g_m^{-1}}B_{g_m^{-1}g_{
m-1}^{-1}}\cdots B_{g_m^{-1}\cdots g_1^{-1}}\,,$$
then 
$$\eqalign{\beta_p(_A\langle x,y\rangle )&=\alpha_{g_1}\alpha_{g_
1}^{-1}\cdots\alpha_{g_m}\alpha_{g_m}^{-1}\alpha_s(_A\langle x,y\rangle 
)\cr
&={}_A\langle\phi_{g_1}\phi_{g_1}^{-1}\cdots\phi_{g_m}\phi_{g_m}^{
-1}\phi_s(x),\phi_{g_1}\phi_{g_1}^{-1}\cdots\phi_{g_m}\phi_{g_m}^{
-1}\phi_s(y)\rangle\cr
&={}_A\langle\phi_p(x),\phi_p(y)\rangle\,.\cr}
$$
Similar calculations show that
Definition \cite{impiso}(a) is also satisfied, which is enough by Lemma 
\cite{relaxlem}.  \eop

\uj Starting with a twisted partial action $(A,G,\alpha ,u)$, Exel [Ex3] builds a 
semidirect product $C^{*}$-algebraic bundle $B$ over $G$ in the sense of Fell.  
He defines [Ex3, Introduction] the crossed product $A\times_{\alpha 
,u}G$ as the 
enveloping $C^{*}$-algebra of the cross sectional algebra $L^1(B)$.  We show 
that the corresponding Busby-Smith twisted action has an isomorphic 
crossed product:  
 
\prop If the Busby-Smith twisted action $(A,S(G),\beta ,w)$ corresponds to 
the twisted partial action $(A,G,\alpha ,u)$ then the crossed products
$A\times_{\alpha ,u}G$ and $A\times_{\beta ,w}S(G)$ are isomorphic.

\proof We are going to show that the Banach $*$-algebras $L_{\beta}
/I_{\beta}$ and 
$L^1(B)$ are isomorphic, which suffices since the crossed products are 
the enveloping $C^{*}$-algebras.  The formula 
$$\phi'(a\chi_{[g_1]\cdots [g_n^{-1}][s]}):=a\chi_s$$
defines a bounded $*$-homomorphism $\phi':L_{\beta}\to L^1(B)$. Since 
$$\phi'(a\chi_{[g_1]\cdots [g_n^{-1}][e]}-a\chi_{[e]})=a\chi_e-a\chi_
e=0\,,$$
$\phi'$ takes $I_{\beta}$ to $0$ and hence determines a bounded $
*$-homomorphism 
$\phi :L_{\beta}/I_{\beta}\to L^1(B)$. Going the other way, the formula 
$$\psi (a\chi_s):=a\chi_{[s]}+I_{\beta}$$
defines a bounded $*$-homomorphism $\psi :L^1(B)\to L_{\beta}$. It is clear that $
\psi\circ\phi$ 
is the identity map. To show that $\psi\circ\phi$ is also the identity map, 
consider $\psi\circ\phi (a\chi_{[g_1]\cdots [g_n^{-1}][s]}+I_{\beta}
)=a\chi_{[s]}+I_{\beta}$. We can choose 
elements $b,c\in A_{[g_1]\cdots [g_n^{-1}][s]}$ such that $a=bc$. Hence 
$$a\chi_{[g_1]\cdots [g_n^{-1}][s]}-a\chi_{[s]}=(b\chi_{[g_1]\cdots 
[g_n^{-1}]}-b\chi_{[e]})*c\chi_{[s]}\in I_{\beta}\,.$$
\Eop

\uj
Using Theorems \cite{moreqbsxp} and \cite{moreqcor} we now have:

\cor Morita equivalent twisted partial actions have Morita equivalent 
crossed products. 

\uj\ 
We now develop the basic theory of covariant representations for 
twisted partial actions.

\mark{tcovrep}
\definition A {\it covariant representation\/} of a twisted partial action 
$(A,G,\alpha ,u)$ is a triple $(\pi ,U,H)$, where $\pi$ is a nondegenerate 
representation of $A$ on the Hilbert space $H$ and for all $s\in 
G$, $U_s$ is a 
partial isometry on $H$ such that
\item{(a)}{$U_s$ has initial space $\pi (A_{s^{-1}})H$ and final space $
\pi (A_s)H$;}
\item{(b)}{$U_sU_t=\pi (u_{s,t})U_{st}$ for all $s,t\in G$;}
\item{(c)}{$\pi (\alpha_s(a))=U_s\pi (a)U_s^{*}$ for all $a\in A_{
s^{-1}}$.}
\uj
Note that we have $U_{s^{*}}=\pi (u_{s^{*},s})U_s^{*}$ for all $s
\in G$. Every covariant 
representation gives a representation of the cross sectional algebra:

\definition The {\it integrated form} $\pi\times U:L_1(B)\to B(H)$ of the covariant 
representation $(\pi ,U)$ is defined by 
$$\big(\pi\times U\big)(x)=\sum_{s\in G}\pi (x(s))U_s\,,$$
where the series converges in norm.
 
\uj The proof of the following proposition is essentially the same as 
that of [Si2, Proposition 3.5].  

\prop$\pi\times U$ is a nondegenerate representation of $L_1(B)$.

\mark{covreplem}
\lemma Let $(A,S(G),\beta ,v)$ be a Busby-Smith twisted action corresponding 
to the twisted partial action $(A,G,\alpha ,u)$.  If $(\pi ,V)$ is a covariant 
representation of $(\beta ,v)$ then $(\pi ,U)$ is a covariant representation of 
$(\alpha ,u)$, where $U_s:=V_{[s]}$ for all $s\in G$.  Conversely, if $
(\pi ,U)$ is a 
covariant representation of $(\alpha ,U)$ then $(\pi ,V)$ is a covariant 
representation of $(\beta ,v)$, where 
$$V_{[g_1][g_1^{-1}]\cdots [g_n][g_n^{-1}][s]}:=P_{g_1}\cdots P_{
g_n}U_s$$
and $P_t$ denotes $\pi (1_{M(A_t)})$ for all $t\in G$.  Moreover this correspondence 
between covariant representations of $(\alpha ,u)$ and $(\beta ,v
)$ is bijective.  

\proof The only nontrivial condition to check for the first part is 
Definition \cite{tcovrep}(b):
$$\eqalign{U_sU_t&=V_{[s]}V_{[t]}=\pi (v_{[s],[t]})V_{[s][s^{-1}]
[st]}=\pi (1_{M(A_{[s][s^{-1}][t]})}u_{s,t})V_{[s][s^{-1}][st]}\cr
&=\pi (u_{s,t})\pi (1_{M(A_{[s][s^{-1}][t]})})V_{[s][s^{-1}][st]}
=\pi (u_{s,t})\pi (v_{[s][s^{-1}],[st]})V_{[s][s^{-1}]}V_{[st]}\cr
&=\pi (u_{s,t})V_{[s][s^{-1}]}V_{[st]}=\pi (u_{s,t})U_{st}\,.\cr}
$$
To show the second part first notice that the $P_t$'s commute since the 
$1_{M(A_t)}$'s are central projections in the double dual of $A$. Therefore 
$V_{[g_1]\cdots [g_n^{-1}][s]}$ is well defined since $P_{g_1}\cdots 
P_{g_n}$ does not depend on the 
order of the idempotents $[g_1][g_1^{-1}],\dots,[g_n][g_n^{-1}]$. It is clear that 
$V_{[g_1]\cdots [g_n^{-1}][s]}$ is a partial isometry. This partial isometry has the 
required final space since
$$\eqalign{\pi (V_{[g_1]\cdots [g_n^{-1}][s]})H&=P_{g_1}\cdots P_{
g_n}U_sH=P_{g_1}\cdots P_{g_n}\pi (A_s)H\cr
&=P_{g_1}\cdots P_{g_n}P_sH=\pi (A_{g_1}\cdots A_{g_n}A_s)H\cr
&=\pi (A_{[g_1]\cdots [g_n^{-1}][s]})H\,.\cr}
$$
We can show that it also has the required initial space by taking 
conjugates.  To check multiplicativity, let $p=[g_1]\cdots [g_m^{
-1}][s]$ and 
$q=[h_1]\cdots [h_n^{-1}][t]$.  Then we have 
$$\eqalign{V_pV_q&=P_{g_1}\cdots P_{g_m}U_sP_{h_1}\cdots P_{h_n}U_
t\cr}
\,.$$
We first simplify a piece of this expression:
$$\eqalign{U_sP_{h_1}&=U_sU_s^{*}U_sU_{h_1}U_{h_1}^{*}\cr
&=P_s\pi (u_{s,h_1})U_{sh_1}U_{h_1}^{*}\cr
&=P_s\pi (u_{s,h_1})P_{sh_1}U_{sh_1}U_{h_1^{-1}}\pi (u_{h_1,h_1^{
-1}})^{*}\cr
&=P_sP_{sh_1}\pi (u_{s,h_1})\pi (u_{sh_1,h_1^{-1}})U_s\pi (u_{h_1
,h_1^{-1}})^{*}\cr
&=P_sP_{sh_1}\pi (u_{s,h_1})U_sU_s^{*}\pi (u_{sh_1,h_1^{-1}})U_s\pi 
(u_{h_1,h_1^{-1}})^{*}\cr
&=\lim_{\lambda}P_sP_{sh_1}\pi (u_{s,h_1})U_s\pi (\alpha_s^{-1}(e_{
\lambda}u_{sh_1,h_1^{-1}}))\pi (u_{h_1,h_1^{-1}})^{*}\,,\cr
&\qquad\hbox{\rm where $e_{\lambda}$ is an approximate identity for $
A_sA_{sh_1}$}\cr
&=\lim_{\lambda}P_sP_{sh_1}\pi (u_{s,h_1})U_s\pi (u_{s^{-1},s}^{*}
\alpha_{s^{-1}}(e_{\lambda}u_{sh_1,h_1^{-1}})u_{s^{-1},s})\pi (u_{
h_1,h_1^{-1}})^{*}\cr
&=\lim_{\lambda}P_sP_{sh_1}\pi (u_{s,h_1})U_s\pi (u_{s^{-1},s}^{*}
\alpha_{s^{-1}}(e_{\lambda})u_{s^{-1},sh_1}u_{h_1,h_1^{-1}})\pi (
u_{h_1,h_1^{-1}})^{*}\cr
&=P_sP_{sh_1}\pi (u_{s,h_1})U_s\pi (u_{s^{-1},s}^{*}u_{s^{-1},sh_
1})\cr
&=\lim_{\mu}P_sP_{sh_1}U_sU_s^{*}\pi (e_{\mu}u_{s,h_1})U_s\pi (u_{
s^{-1},s}^{*}u_{s^{-1},sh_1})\cr
&=\lim_{\mu}P_sP_{sh_1}U_s\pi (\alpha_s^{-1}(e_{\mu}u_{s,h_1}))\pi 
(u_{s^{-1},s}^{*}u_{s^{-1},sh_1})\cr
&=\lim_{\mu}P_sP_{sh_1}U_s\pi (u_{s^{-1},s}^{*}\alpha_{s^{-1}}(e_{
\mu}u_{s,h_1})u_{s^{-1},s})\pi (u_{s^{-1},s}^{*}u_{s^{-1},sh_1})\cr
&=\lim_{\mu}P_sP_{sh_1}U_s\pi (u_{s^{-1},s}^{*}\alpha_{s^{-1}}(e_{
\mu})u_{s^{-1},s}u_{e,h_1}u_{s^{-1},sh_1}^{*})\pi (u_{s^{-1},sh_1}
)\cr
&=P_sP_{sh_1}U_s\,.\cr}
$$
Repeating this calculation $n-1$ times we have 
$$\eqalign{V_pV_q&=P_{g_1}\cdots P_{g_m}P_sP_{sh_1}P_{sh_2}\cdots 
P_{sh_n}U_sU_t\cr
&=P_{g_1}\cdots P_{g_m}P_sP_{sh_1}\cdots P_{sh_m}\pi (u_{s,t})U_{
st}\cr
&=\pi (v_{p,q})V_{[g_1]\cdots [g_m^{-1}][s][s^{-1}][sh_1]\cdots [
sh_m^{-1}][st]}\cr
&=\pi (v_{p,q})V_{pq}\,.\cr}
$$
Finally we check the covariance condition. If $p=[g_1]\cdots [g_m^{
-1}][s]$ and 
$a\in A_{p^{*}}$ then 
$$\eqalign{\pi (\beta_p(a))&=\pi (\alpha_s(a))\cr
&=\pi (\alpha_{g_1}\alpha_{g_1}^{-1}\cdots\alpha_{g_n}\alpha_{g_n}^{
-1}\alpha_s(a))\cr
&=U_{g_1}\pi (u_{g_1^{-1},g_1}^{*}\alpha_{g_1^{-1}}\cdots\alpha_{
g_n}\alpha_{g_n}^{-1}\alpha_s(a)u_{g_1^{-1},g_1})U_{g_1}^{*}\cr
&=U_{g_1}\pi (u_{g_1^{-1},g_1}^{*})U_{g_1^{-1}}\pi (\alpha_{g_2}\cdots
\alpha_{g_n}\alpha_{g_n}^{-1}\alpha_s(a))U_{g_1^{-1}}^{*}\pi (u_{
g_1^{-1},g_1})U_{g_1}^{*}\cr
&=U_{g_1}U_{g_1}^{*}\pi (\alpha_{g_2}\cdots\alpha_{g_n}\alpha_{g_
n}^{-1}\alpha_s(a))U_{g_1^{-1}}^{*}U_{g_1^{-1}}\cr
&=\cdots\cr
&=P_{g_1}\cdots P_{g_n}U_s\pi (a)U_s^{*}P_{g_n}^{*}\cdots P_{g_1}^{
*}\cr
&=V_p\pi (a)V_p^{*}\,.\cr}
$$
It is clear from the construction that the correspondence is bijective.
\eop

\prop If $(A,G,\alpha ,u)$ is a twisted partial action then $(\pi 
,U)\mapsto\pi\times U$ is a 
bijective correspondence between covariant representations of $(\alpha 
,u)$ 
and nondegenerate representations of the crossed product $A\times_{
\alpha ,u}G$.  

\proof We know that there is an isomorphism $\phi$ between $A\times_{
\beta ,v}S(G)$ 
and $A\times_{\alpha ,u}G$ where $(A,S(G),\beta ,v)$ is the corresponding semigroup 
action.  We also know that there is a bijective correspondence 
$\Psi\mapsto (\pi^{\Psi},V^{\Psi})$ between nondegenerate representations of $
A\times_{\beta ,v}S(G)$ and 
covariant representations of $(\beta ,v)$ such that $\Psi =\pi^{\Psi}
\times V^{\Psi}$.
We define a bijective correspondence $\Phi\mapsto (\pi^{\Phi},U^{
\Phi})$ 
between nondegenerate representations of $A\times_{\alpha ,u}G$ and covariant 
representations of $(\alpha ,u)$ satisfying $\Phi =\pi^{\Phi}\times 
U^{\Phi}$ using the following 
diagram:
$$\matrix{\hbox{\rm Rep}\,(A\times_{\alpha ,u}G)&\leftrightarrow&\hbox{\rm Rep}\,
(A\times_{\beta ,v}S(G))&\qquad\qquad&\Phi&\leftrightarrow&\Psi =
\Phi\circ\phi\cr
\cr
&&\updownarrow&&&&\updownarrow\cr
\cr
\hbox{\rm CovRep}\,(\alpha ,u)&\leftrightarrow&\hbox{\rm CovRep}\,
(\beta ,v)&&(\pi^{\Phi},U^{\Phi})&\leftrightarrow&(\pi^{\Psi},V^{
\Psi})\cr}
$$
If $\Phi$ is a nondegenerate representation of $A\times_{\alpha ,
u}G$ then $\Psi =\Phi\circ\phi$ is a 
nondegenerate representation of $A\times_{\beta ,v}S(G)$ and so $
\Psi =\pi^{\Psi}\times V^{\Psi}$.  Let 
$(\pi^{\Phi},U^{\Phi})$ be the covariant representation of $(\alpha 
,u)$ corresponding to 
$(\pi^{\Psi},V^{\Psi})$ as in Lemma \cite{covreplem}.  If $a\in A_
s$ then 
$$\pi^{\Phi}\times U^{\Phi}(a\delta_s)=\pi^{\Phi}(a)U^{\Phi}_s=\pi^{
\Psi}(a)V^{\Psi}_{[s]}=\Psi (a\delta_{[s]})=\Phi (a\delta_s)$$
and so $\pi^{\Phi}\times U^{\Phi}=\Phi$.
\eop

\newsection{Connection with crossed products by Hilbert bimodules}

\uj Recall from [AEE] that the crossed product $A\times_{\alpha}{\bf Z}$ of 
the partial action $(A,{\bf Z},\alpha )$ is isomorphic to the crossed 
product $A\times_X{\bf Z}$ of $A$ by the Hilbert bimodule $_AX_A$, where $
X$ 
is the vector space $A_1$ with module structure 
$$a\cdot j:=aj,\quad j\cdot a:=\alpha_1(\alpha_1^{-1}(j)a)$$
and inner products
$$_A\langle j,k\rangle :=jk^{*},\quad\langle j,k\rangle_A:=\alpha_
1^{-1}(j^{*}k)$$
for $j,k\in A_1$ and $a\in A$. In other words, we can get $_AX_A$ by converting 
the standard $A_1-A_1$ imprimitivity bimodule $A_1$ into an $A_1-
A_{-1}$ 
imprimitivity bimodule via the isomorphism $\alpha_1$, then extending it 
canonically to a Hilbert $A-A$ bimodule.

\definition The Hilbert bimodules $_AX_A$ and $_BY_B$ are called {\it Morita }
{\it equivalent\/} if there is an isomorphism $(\hbox{\rm id},\phi 
,\hbox{\rm id})$ between the Hilbert 
bimodules $X\otimes_AM$ and $M\otimes_BY$ for some imprimitivity bimodule $_
AM_B$.  

\uj Abadie, Eilers and Exel show that if $_AX_A$ and $_B$$Y_B$ are Morita 
equivalent bimodules then the crossed products $A\times_X{\bf Z}$ and $
B\times_Y{\bf Z}$ are 
Morita equivalent.  They note that Hilbert bimodules corresponding to 
Morita equivalent actions of ${\bf Z}$ are Morita equivalent.  We show that 
the Morita equivalence of Hilbert bimodules corresponding to partial 
actions of ${\bf Z}$ is equivalent to the Morita equivalence of the 
partial actions, in the sense of Definition \cite{moreqpar}.  

Suppose we have two partial actions $(A,\alpha ,{\bf Z})$ and $(B
,\beta ,{\bf Z})$ with 
corresponding Hilbert bimodules $_AX_A$ and $_BY_B$.  We show that the 
two notions of Morita equivalence of the actions coincide.  

\prop The partial actions $(A,\alpha ,{\bf Z})$ and $(B,\beta ,{\bf Z}
)$ are Morita 
equivalent if and only if the corresponding Hilbert 
bimodules $_A$$X_A$ and $_BY_B$ are Morita equivalent.  

\proof If $_AM_B$ is an imprimitivity bimodule then 
$$\eqalign{\overline {\hbox{\rm span}}\,\langle M\otimes_BY,M\otimes_
BY\rangle_B&=\overline {\hbox{\rm span}}\,\langle Y,\langle M,M\rangle_
B\cdot Y\rangle_B\cr
&=\overline {\hbox{\rm span}}\,\beta_1^{-1}(B_1^{*}\langle M,M\rangle_
BB_1)=B_{-1}\,,\cr}
$$
hence the imprimitivity bimodule corresponding to $M\otimes_BY$ is of the 
form $_D(M\otimes_BY)_{B_{-1}}$ for some closed ideal $D$ of $A$.  Similarly, the 
imprimitivity bimodule corresponding to $_A(X\otimes_AM)_B$ is of the form 
$_{A_1}(X\otimes_AM)_C$ for some closed ideal $C$ of $B$.  It is routine to check 
that the map $m\otimes l\mapsto m\cdot l$ for $m\in M$ and $l\in 
B_1$ extends to a map 
$\nu :M\otimes_BY\to M\cdot B_1$ such that $(\hbox{\rm id},\nu ,\beta_
1)$ is an isomorphism between 
$_D(M\otimes_BY)_{B_{-1}}$ and the imprimitivity subbimodule $_D(
M\cdot B_1)_{B_1}$ of $_AM_B$.  
Similarly, the map $j\otimes m\mapsto\alpha_{-1}(j)\cdot m$ for $
j\in A_1$ and $m\in M$ extends to a 
map $\mu :X\otimes_AM\to A_{-1}\cdot M$ such that $(\alpha_{-1},\mu 
,\hbox{\rm id})$ is an isomorphism 
between $_{A_1}(X\otimes_AM)_C$ and the imprimitivity subbimodule $_{
A_{-1}}(A_{-1}\cdot M)_C$ 
of $_AM_B$.  

Suppose now that the Hilbert bimodules $X$ and $Y$ are Morita 
equivalent.  Then by Lemma \cite{idfiidlem} and the above there 
exists an imprimitivity bimodule $_AM_B$ and an isomorphism $(\hbox{\rm id}
,\psi ,\hbox{\rm id})$ 
between the imprimitivity bimodules $_{A_1}(X\otimes_AM)_C$ and $_
D(M\otimes_BY)_{B_{-1}}$.  
Then $A_1=D$ and $C=B_{-1}$ and so $(\alpha_1,\nu\circ\psi\circ\mu^{
-1},\beta_1)$ is an isomorphism 
between $_{A_{-1}}(A_{-1}\cdot M)_{B_{-1}}$ and $_{A_1}(M\cdot B_
1)_{B_1}$.  This implies that 
$(A,X,{\bf Z})\sim_{X,\phi}(B,Y,{\bf Z})$, where $\phi_n$$:=(\nu\circ
\psi\circ\mu^{-1})^n$ for $n\in {\bf Z}\setminus \{0\}$.  The 
situation can be visualized by the following diagram:  
\def \sc{\scriptstyle} 
\input paper4.fi1

Going the other way, if $(A,\alpha ,{\bf Z})\sim_{M,\phi}(B,\beta 
,{\bf Z})$ then $\phi_1$ is 
an isomorphism between $_{A_{-1}}($$A_{-1}\cdot M)_{B_{-1}}$ and $_{
A_1}($$M\cdot B_1)_{B_1}$. So $C=B_{-1}$, 
$D=A_1$ and $(\hbox{\rm id},\nu^{-1}\circ\phi_1\circ\mu ,\hbox{\rm id}
)$ is an isomorphism between     
$_{A_1}(X\otimes_AM)_C$ and $_D(M\otimes_BY)_{B_{-1}}$ and so the Hilbert bimodules $
X$ and $Y$ 
are Morita equivalent by Lemma \cite{idfiidlem}.
\eop

\bigskip\goodbreak

\centerline{\nagy References}
\nobreak
\medskip
\nobreak

\def \i#1#2{\itemitem {\hbox to .45in{#1\hfil}} {#2} }
\rm

\i{[AEE]}{B. Abadie, S. Eilers and R. Exel, {\it Morita 
equivalence for crossed products by Hilbert $C^*$-bimodules},
Trans. Amer. Math. Soc. {\bf 350}
(1998), no.~8, 3043--3054.}
\i{[BGR]}{L. Brown, P. Green and M. Rieffel,
{\it Stable isomorphism and strong Morita equivalence of 
$C^{*}$-algebras}, Pac. J. Math. {\bf 71}(2) (1977), 349--363.}
\i{[Com]}{F. Combes, {\it Crossed products and Morita 
equivalence}, Proc. London. Math. Soc. {\bf 49} (1984), 
289--306.}
\i{[CMW]}{R. Curto, P. Muhly and D. Williams, 
{\it Cross products of strongly Morita equivalent $C^{*}$-algebras},
Proc. Amer. Math. Soc. {\bf 90}(4) (1984), 528--530. }
\i{[Ech]}{S. Echterhoff, 
{\it Morita equivalent twisted actions and a new version of the 
Packer-Raeburn stabilization trick}, 
J. London Math. Soc. (2) {\bf 50} (1994), 170--186.}
\i{[ER]}{S. Echterhoff and I. Raeburn, {\it Multipliers of 
imprimitivity bimodules and Morita equivalence of 
crossed products}, Math. Scand. {\bf 76} (1995), 289--309.}
\i{[Ex1]}{R. Exel, {\it Circle actions on $C^{*}$-algebras, 
partial automorphisms
and a generalized Pimsner-Voiculescu exact sequence}, J.  
Funct.  Anal. {\bf 122}  (1994), 361--401.}
\i{[Ex2]}{R. Exel, {\it Twisted partial actions: a classification 
of regular $C^{*}$-algebraic bundles},
Proc. London Math. Soc. (3) {\bf 74} (1997), no.~2, 417--443.}
\i{[Ex3]}{R. Exel, \it Partial actions of groups and actions 
of inverse semigroups\rm, Proc. Amer. Math. Soc.
{\bf 126} (1998), no.~12, 3481--3494.}
\i{[JT]}{K. Jensen and K. Thomsen, {\it Elements of KK-theory}, 
Birkh\"auser, Boston, 1991. }
\i{[Kal]}{S. Kaliszewski, {\it Morita equivalence methods 
for twisted $C^*$-dynamical systems}, Ph.D. thesis, Dartmouth College, 
1994.} 
\i{[Lan]}{E. Lance, {\it Hilbert $C^*$-modules, a toolkit for operator
algebraists}, Cambridge University Press, Cambridge, 1995.}
\i{[Pat]}{A. Paterson, {\it r-Discrete $C^*$-algebras as covariant 
$C^*$-algebras}, Groupoid Fest lecture notes, Reno, 1996.}
\i{[Rie]}{M. Rieffel, {\it Unitary representations of group 
extensions: an algebraic approach to the theory of 
Mackey and Blattner}, Advances in Mathematics 
Supplementary Studies, {\bf 4} (1979), 43--81.}
\i{[Si1]}{N. Sieben, $C^{*}${\it-crossed products by partial 
actions and actions of inverse semigroups},
J. Operator Theory {\bf 39} (1998), no.~2, 361--393.}
\i{[Si2]}{N. Sieben, {\it $C^{*}$-crossed products by twisted 
inverse semigroup actions},
J. Austral. Math. Soc. Ser. A {\bf 63} (1997), no.~1, 32--46.}

\bye

%% file: makro.tex
\tolerance=10000
\magnification=\magstephalf
\def\nagy{\font\caps=cmcsc10\caps}
\def\kis{\font\kisf=cmr5\kisf}

\def\sethoff{\hoffset=0.5 true in}
\sethoff  
\hsize=6.0 true in
\voffset=.3 true in
\vsize=8.5 true in  

\newif\ifduplexpr
\duplexprfalse

\def\duplex{\hoffset=0.25 true in \duplexprtrue \marginfalse}


\def\uj {\bigskip \rm}


\def\redef#1#2{\expandafter\ifx\csname #1\endcsname\relax
    \expandafter\edef\csname #1\endcsname{#2} 
    \else \message{redefinition of '\string#1'
    } \fi }

\newif\ifexist

\def\testfile#1{
\openin 0=#1
\ifeof 0
\existfalse
\else
\existtrue
\fi
\closein 0
}


\testfile{cite.inc}
\ifexist
\input cite.inc
\else
\message{!!!!!!!!!! One more pass needed for references !!!!!!!!!}
\fi


\immediate\openout 0=cite.inc


\immediate\openout 2=conten.inc
\def\writeitem#1#2{\write2
  {\string
  \line{#1\string\nagy{}#2\string\rm\string\leaderfill\string\quad\folio}}}


\def\boxit#1{\vbox{\hrule\hbox{\vrule#1\vrule}\hrule}}


\newif\ifmargin
\margintrue


\newcount\sorszam \sorszam=0
\def\sorszaminc {\advance\sorszam by1 
\ifnum \section=0 \else \the\section.\fi
\the\sorszam. }

\def\authorstr{}
\def\shorttitlestr{}

\def\author#1{\def\authorstr{#1} 
 \bigskip
 \centerline{#1}
}

\def\abstract#1{\bigskip
 \centerline{\nagy Abstract}\medskip
 \vbox{\narrower\narrower\noindent#1}
}

\def\shorttitle#1{\def\shorttitlestr{#1}}

\def\oldalszam{
\headline={
 \nagy
 \ifnum \folio > 1
 \ifodd \folio {\hfil \shorttitlestr \hfil\folio} 
 \else
 {\folio\hfil \authorstr \hfil}
 \fi
 \else \hfil
 \fi
}
}

\nopagenumbers
\oldalszam


\newcount\section \section=0
\def\newsection#1{
      \goodbreak
      \advance\section by1 \sorszam=0 
      \medskip\bigskip\centerline{\nagy \the\section. #1}
      \nobreak\nobreak
     \writeitem{\the\section.  }{#1}
     }
\newbox\labelboxx
\def\ugor{\goodbreak
\bigskip
\ifmargin
\hskip-1in {\box\labelboxx}
\vskip-\baselineskip
\fi
\noindent\bf}
\def\skippy{\enskip}

\def\definition { \ugor Definition \sorszaminc \rm}   
\def\theorem {    \ugor Theorem \sorszaminc \sl}   
\def\prop {       \ugor Proposition \sorszaminc \sl}
\def\lemma {      \ugor Lemma \sorszaminc \sl}      
\def\cor  {       \ugor Corollary \sorszaminc \sl} 
      
\def\proof { \medskip {\noindent \it Proof.\skippy}\rm}

\def\nullbox{\setbox0=\null \ht0=5pt \wd0=5pt \dp0=0pt \box0}
\def\eop {\hfill\hbox{\relax}\hfill \boxit{\nullbox}\goodbreak}    
\def\Eop {\vskip -\baselineskip \vskip -\belowdisplayshortskip \eop}


\def\label#1{  
\ifmargin
\vskip0in\hskip-1in {\kis#1}
\vskip-\baselineskip
\fi
\immediate
\write0{\string\redef{\string#1}{{\string\rm\the\section.\the\sorszam}}}}

\def\mark#1{  
\advance\sorszam by 1
\setbox\labelboxx=\hbox{\kis #1}
\immediate
\write0{\string\redef{\string#1}{{\string\rm
\ifnum \section=0 \else\the\section.\fi
\the\sorszam}}}
\advance\sorszam by -1
}

\def\labelref#1{ 
\immediate
\write0{\string\redef{\string#1}{{\string\rm\the\sorszam}}}}

\def\cite#1{\expandafter\ifx\csname#1\endcsname\relax
    ???
    \message{!!!!!!! Missing reference '\string#1' !!!!!!!}
    \else
\csname#1\endcsname
\fi}


\def\title#1{\shorttitle{#1}\centerline{\nagy #1}}


\def\i#1#2#3#4{\itemitem{\hbox to .5in{#1\hfil}}{#2, {\it #3}, #4}}



%% file: cite.inc
\redef{fullemma}{{\rm2.1}}
\redef{impiso}{{\rm2.4}}
\redef{relaxlem}{{\rm2.5}}
\redef{idfiidlem}{{\rm2.6}}
\redef{morequ}{{\rm3.2}}
\redef{fiax}{{\rm3.3}}
\redef{linkaction}{{\rm3.5}}
\redef{moreqbsxp}{{\rm3.8}}
\redef{unequ}{{\rm3.9}}
\redef{partialaction}{{\rm4.1}}
\redef{moreqpar}{{\rm4.2}}
\redef{moreqcor}{{\rm4.3}}
\redef{tcovrep}{{\rm4.6}}
\redef{covreplem}{{\rm4.9}}

%% file: smalmtrx.tex
\def\smallmatrix#1{\null\,\vcenter{\baselineskip=0.7\baselineskip
    \ialign{\hfil$\scriptstyle##$\hfil&&$\,\,$\hfil$\scriptstyle##$\hfil\crcr
      \mathstrut\crcr\noalign{\kern-\baselineskip}
      #1\crcr\mathstrut\crcr\noalign{\kern-\baselineskip}}}\,
      \normalbaselines}

%% file: paper4.fi1
$$\matrix{_A(X\otimes_AM)_B&&_{A_1}(X\otimes_AM)_C&\buildrel{\sc(\alpha_{-1}
,\mu ,\hbox{\fiverm id})}\over{\relbar\joinrel\relbar\joinrel\relbar\joinrel\relbar\joinrel\longrightarrow}&_{A_{-1}}(A_{-1}\cdot M)_C\cr
\cr
\Big\downarrow&&\quad\Big\downarrow\sc(\hbox{\rm \fiverm id}\,,\psi 
,\hbox{\rm \fiverm id}\,)&&\qquad\quad\Big\downarrow\sc(\alpha_1,
\phi_1,\beta_1)\cr
\cr
_A(M\otimes_BY)_B&&_D(M\otimes_BY)_{B_{-1}}&\buildrel{\sc(\hbox{\fiverm id},\nu 
,\beta_1)}\over{\relbar\joinrel\relbar\joinrel\relbar\joinrel\relbar\joinrel\longrightarrow}&_D(M\cdot B_1)_{B_1}\cr}
$$